\documentclass[1p]{elsarticle}
\usepackage{amsmath}
\usepackage{amssymb}
\usepackage{mathtools}
\usepackage{bm}
\usepackage{graphicx,float}
\usepackage{tikz}
\usepackage[title]{appendix}
\usepackage{hyperref}

\makeatletter
\newcommand*\bigcdot{\mathpalette\bigcdot@{.6}}
\newcommand*\bigcdot@[2]{\mathbin{\vcenter{\hbox{\scalebox{#2}{$\m@th#1\bullet$}}}}}
\makeatother

\renewcommand{\vec}{\boldsymbol}
\newcommand{\del}{\vec{\nabla}}
\newcommand{\norm}[1]{\left\lVert#1\right\rVert}
\newcommand{\dd}{\,\text{d}}

\newcommand{\dS}{\partial_{S}}
\newcommand{\dPhi}{\partial_{\Phi}}

\newcommand{\dZ}{\partial_{Z}}

\newcommand{\eS}{\vec{\hat{e}}_{S}}
\newcommand{\ePhi}{\vec{\hat{e}}_{\Phi}}
\newcommand{\eZ}{\vec{\hat{e}}_{Z}}

\newcommand{\ds}{\partial_{s}}
\newcommand{\dphi}{\partial_{\phi}}
\newcommand{\deta}{\partial_{\eta}}
\newcommand{\dzeta}{\partial_{\zeta}}
\newcommand{\htilde}{\, \widetilde{h}}
\newcommand{\stilde}{\, \widetilde{s}}

\usepackage{xcolor}

\makeatletter
\renewcommand*\env@matrix[1][\arraystretch]{%
  \edef\arraystretch{#1}%
  \hskip -\arraycolsep
  \let\@ifnextchar\new@ifnextchar
  \array{*\c@MaxMatrixCols c}}
\makeatother

\begin{document}

\begin{frontmatter}

\title{Gyroscopic polynomials}

\author[cuboulder]{Abram C. Ellison}
\author[cuboulder]{Keith Julien}

\affiliation[cuboulder]{organization={University of Colorado Boulder Department of Applied Mathematics},
                        city={Boulder},
                        state={CO 80309},
                        country={USA}}

\begin{abstract}
Gyroscopic alignment of a fluid occurs when flow structures align with the rotation axis.
This often gives rise to highly spatially anisotropic columnar structures that in combination with complex domain
boundaries pose challenges for efficient numerical discretizations and computations.
We define gyroscopic polynomials to be three-dimensional polynomials expressed in a coordinate system that conforms to rotational alignment.
We remap the original domain with radius-dependent boundaries onto a right cylindrical or annular domain to create the computational domain in this coordinate system.
We find the volume element expressed in gyroscopic coordinates leads naturally to a hierarchy of orthonormal bases.
We build the bases out of Jacobi polynomials in the vertical and generalized Jacobi polynomials in the radial.  Because these coordinates explicitly
conform to flow structures found in rapidly rotating systems the bases represent fields with a relatively small number of modes.
We develop the operator structure for one-dimensional semi-classical orthogonal polynomials as a building block for differential operators in the full
three-dimensional cylindrical and annular domains.  The differential operators of generalized Jacobi polynomials generate a sparse
linear system for discretization of differential operators acting on the gyroscopic bases.  This enables efficient
simulation of systems with strong gyroscopic alignment.
\end{abstract}



\end{frontmatter}

\section{Introduction}
Gyroscopic alignment of a fluid occurs when flow structures align with the rotation axis.
It is typical that such flows are highly spatially anisotropic columnar structures irrespective of the geometry that confines them.
Natural questions that arise in the presence of spatial anisotropy and complex geometry are (i) what functional bases
are optimal for accurate representation and computation and (ii) can discretization of differential operators remain
sparse such that fast algorithmic strategies are possible.
We define a \emph{gyroscopic polynomial} to be a three-dimensional polynomial expressed in a coordinate system that conforms to this rotational alignment.
For example, the recently developed Spherinder basis \cite{Ellison_Julien_Vasil_2022} is an orthonormal
gyroscopic polynomial basis for the unit ball.  This basis eschews spherical coordinates to better represent flows found in rapidly rotating fluids.
A typical discretization of the unit ball is an expansion in spherical harmonics:
\begin{equation}
f(r,\theta,\phi) = \sum_{m,l} Y_{l,m}(\theta,\phi) F_{l,m}(r).
\end{equation}
Though the coordinate system $(r,\theta,\phi)$ is orthogonal, gyroscopically aligned field contours parallel to
$\vec{\hat{e}}_{Z} = \cos \theta \vec{\hat{e}}_{r} - \sin \theta \vec{\hat{e}}_{\phi}$
are functions of both the radius $r$ and colatitude $\theta$.
This means axially aligned fields require a large number of modes before the expansion converges to a sufficient accuracy.
The Spherinder methodology \cite{Ellison_Julien_Vasil_2022} demonstrates how to design spectral bases for gyroscopic coordinates to
get around this issue.  Not only do the bases represent gyroscopic fields with relatively few modes, but differential operators acting on
these basis functions have a sparse matrix representation.  Spherinder basis discretization is therefore fast and efficient for gyroscopic flows in the sphere geometry.
In this paper we extend this methodology to define gyroscopic polynomials on stretched cylinders and annuli with polynomial bounding surface $h(s)$.

Cylindrical and annular geometries are common in laboratory settings for rotating flows \cite{Lonner_Aggarwal_Aurnou_2022}.
Experimenters design an apparatus to emulate instabilities and turbulent motions relevant to geophysical fluid flow regimes.
The domain is typically a cylinder or annulus of fluid under rapid rotation.
Often the upper surface of the fluid is left free and so forms an equipotential surface - a parabolic contour in this case. 
We will define the \emph{gyroscopic coordinate system} that follows these curved upper surfaces.  
We demonstrate the coordinate system leads naturally to an orthogonal basis of polynomials as in the sphere \cite{Ellison_Julien_Vasil_2022}.
These basis functions can describe fields contained in spheres, ellipsoids, parabolic cylinders and more exotic domains. We generalize the bases one step
further to stretched annular domains with inner radius $S_{i}$ and outer radius $S_{o}$.

Jacobi polynomials that frequently arise in the discretization of partial differential equations (PDEs) are fundamental in the application of gyroscopic polynomials.
A Jacobi polynomial $P_{n}^{(a,b)}(z)$ is orthogonal
under the weight function
\begin{equation}
w^{(a,b)}(z) \triangleq (1-z)^{a} (1+z)^{b}
\end{equation}
such that
\begin{equation}
\left\langle P_{n}^{(a,b)}(z), P_{m}^{(a,b)}(z) \right\rangle_{w^{(a,b)}} \triangleq \int_{-1}^{1} \dd z \, w^{(a,b)}(z) P_{n}^{(a,b)}(z) P_{m}^{(a,b)}(z)
\, \propto \, \delta_{m,n},
\end{equation}
where $\delta_{m,n}$ is the Kronecker delta symbol.  The Legendre polynomials correspond to $(a,b) = (0,0)$ and the Chebyshev polynomials
have $(a,b) = (-1/2,-1/2)$.  A key property of Jacobi polynomials is they form a closed set under differentiation, since
\begin{equation}
\frac{d}{d z} P_{n}^{(a,b)}(z) \propto P_{n-1}^{(a+1,b+1)}(z).
\end{equation}
This sparse derivative representation of $P_{n}^{(a,b)}$ in the $P_{m}^{(a+1,b+1)}$ basis is the catalyst for efficient numerical methods.
In typical problems one uses the volume element of the coordinate system to induce an inner product.  In some cases this volume element
leads directly to a Jacobi-type weight function.  This was the case in \cite{Ellison_Julien_Vasil_2022,Vasil_Burns_Lecoanet_Olver_Brown_Oishi_2016,Vasil_Lecoanet_Burns_Oishi_Brown_2019},
where the authors employed a hierarchy of Jacobi polynomials to discretize cylindrical and spherical coordinate systems.

Classical Jacobi polynomials fall short when generalizing cylindrical domains to have either (i) a polynomial bounding surface or (ii) a finite inner radius $S_{i} > 0$.
Careful study of the volume elements in these geometries reveals the requirement for one or more \emph{augmenting polynomial factors} in the weight function relative
to $w^{(a,b)}$.  To proceed we must extend standard Jacobi polynomials to a more general class of weight functions that accommodates these polynomial
factors, and so we define the \emph{generalized Jacobi polynomials}.  As with Jacobi polynomials, generalized Jacobi polynomials admit sparse derivative operators.
The total degree of the augmenting polynomial factors in the weight function determines the operator bandwidth, so though not as sparse as the classical operators,
the discretized matrix system still has only $\mathcal{O}((d+1) \times N)$ terms, where $d$ is the total augmenting degree and $N$ is the maximum polynomial degree of the truncated system.

Let $(S,\Phi,Z)$ be the standard cylindrical coordinate system and let $h(S)$ denote the polynomial height function of the domain of interest.
When we seek a gyroscopic basis we mean in particular we are looking for approximating functions of the form
\begin{equation}
\Psi_{m,l,k}(S, \Phi, \eta) = e^{i m \Phi} P_{l}\left(\eta\right) Q^{(m,l)}_{k}(S),
\end{equation}
where we made the substitution $Z = \eta \, h(S)$ for the vertical coordinate.
Here $P_{l}$ is a yet-determined degree $l$ polynomial while $Q_{k}^{(m,l)}$ is a radial function parameterized by the azimuthal wavenumber $m$ and vertical degree $l$.
We would like to represent fields in the domain with linear combinations of the basis functions:
\begin{equation}
f(S, \Phi, \eta) = \sum_{m = -\infty}^{\infty} \sum_{l=0}^{\infty} \sum_{k=0}^{\infty} \Psi_{m,l,k}(S,\Phi,\eta) \widehat{F}_{m,l,k}.
\end{equation}
Since gyroscopic alignment suppresses vertical complexity, an aligned field represented in these coordinates will have a rapidly converging expansion in the vertical $l$ index.
We may therefore truncate $l$ to some relatively small $L_{\text{max}}$ compared with the radial truncation length $N_{\text{max}}$ when we discretize the field.
This is the first major motivating factor behind the gyroscopic coordinate system and cannot be overemphasized: \emph{gyroscopically aligned fields admit an expansion in relatively few gyroscopic modes}.

The next step is to identify appropriate classes of functions for $P_{l}(\eta)$ and $Q^{(m,l)}_{k}(S)$.
The $(S, \Phi, \eta)$ coordinates system explicitly separates dynamics in the stretched vertical and radial directions.
When expressing the volume element of this coordinate system we uncover the classical
Jacobi polynomial weight function for the $\eta$ coordinate so that $P_{l}(\eta) \triangleq P_{l}^{(\alpha,\alpha)}(\eta)$ for some $\alpha$.
The radial dependence is more subtle.  First note that the functions $Q^{(m,l)}_{k}(S)$ are \emph{parameterized by $m$ and $l$}.  This parameterization
allows us to (i) guarantee regularity throughout the domain in the presence of coordinate singularities and (ii) ensure the basis $\Psi$ can be expressed as a Cartesian $(x,y,z)$ polynomial.  
We will show that generalized Jacobi polynomials multiplied by the function $S^{m} \, h(S)^{l}$ satisfy the above two conditions.
The three-dimensional basis will inherit the sparsity of differential operators of its one-dimensional polynomial building blocks.
This is the second motivating factor for the gyroscopic basis: \emph{the gyroscopic bases decompose into
one-dimensional (possibly generalized) Jacobi polynomials and hence have sparse matrix representations of algebraic and differential operators}.

In Section~\ref{sec:stretched_cylindrical_coordinates} we describe the gyroscopically aligned coordinate system for cylindrical geometries
and in Section~\ref{sec:stretched_annular_coordinates} we describe the generalization to annular domains.
We develop the generalized Jacobi polynomials in Section~\ref{sec:generalized_jacobi} to be used as building blocks for the gyroscopic polynomials.
In Section~\ref{sec:basis} we define the \emph{gyroscopic bases}, a hierarchy of basis functions used to represent scalar and vector fields.  Section~\ref{sec:discretization}
details how the bases lead to sparse matrix operators for all differential operators needed in fluid dynamics.
We put the bases to the test in Section~\ref{sec:eigenproblems} where we solve the damped inertial waves eigenvalue problem in various geometric
configurations.  We wrap up the paper in Section~\ref{sec:conclusion}.

\section{The Stretched Coordinate Systems}\label{sec:stretched_coordinates}
We adopt the stretched coordinate system ${ \left( s, \phi, \eta \right) }$ related to Cartesian ${(x, y, z)}$ coordinates by
\begin{equation}
\begin{aligned}
x &= s \cos \phi, \\
y &= s \sin \phi, \\
z &= \eta \, h(s)
\end{aligned}
\end{equation}
where
\begin{equation}
\phi \in [0, 2 \pi), \hspace{4ex} \eta \in [-1,1], \hspace{4ex} h(s) \ge 0.
\end{equation}
The domain of $s$ depends on the geometry, where 
\begin{equation}
\begin{aligned}
&s \in \left[ 0, S_{o} \right] \hspace{4ex}& \textrm{(Cylinder)}, \\
&s \in \left[ S_{i}, S_{o} \right] \hspace{4ex}& \textrm{(Annulus)}.
\end{aligned}
\end{equation}
These coordinates have a \emph{polar singularity} along the $z$ axis ($s = 0$).  The coordinate systems also support an \emph{outer equatorial singularity},
sphere-type decay of the height function as $s \to S_{o}$, and the annulus permits an \emph{inner equatorial singularity} with similar decay as $s \to S_{i}$.
Presence of these singularities depends on particular choice of the height function $h(s)$.

Let $(S, \Phi, Z)$ denote standard cylindrical coordinates and ${ \left(\eS, \ePhi, \eZ\right) }$ the associated unit vectors.
Letting $h'(s) = \ds h(s)$, we find partial derivatives transform as
\begin{equation}
\begin{aligned}
\dS &= \ds - \frac{h'(s)}{h(s)} \eta \deta, \\
\dPhi &= \dphi, \\
\dZ &= \frac{1}{h(s)} \deta.
\end{aligned}
\end{equation}
The coupling of the partial derivatives demonstrates that the coordinate vectors $\partial_{s}$ and $\partial_{\eta}$ aren't
orthogonal.  This coupling between stretched coordinates will play out in the sparsity structure of differential operators acting on the gyroscopic bases.
The volume element in our stretched coordinates is
\begin{equation}
\dd V = h(s) s \dd s \dd \phi \dd \eta.
\label{eqn:volume_element_cylinder}
\end{equation}
The key to defining gyroscopic polynomials is using a weighted version of $\dd V$ dictated by $h(s)$ to induce a hierarchy of orthonormal 3D bases.
These bases naturally conform to coordinate singularities and nest in a particular way that yields a sparse linear algebraic structure
for discretized PDEs.

\subsection{Stretched Cylindrical Coordinates}\label{sec:stretched_cylindrical_coordinates}
Due to well-known regularity conditions \cite{Boyd_2001,Matsushima_Marcus_1995,Boyd_Yu_2011,Vasil_Burns_Lecoanet_Olver_Brown_Oishi_2016} in polar coordinates, a scalar field $f$ proportional to $e^{i m \phi}$ decays like 
$s^{m}$ as $s \to 0$, and parity yields
\begin{equation}
f \sim e^{i m \phi} s^{m} F \left( s^{2} \right).
\end{equation}
To explicitly handle the parity condition and to map the physical stretched cylindrical domain onto the computational cylinder we define the variable $t \in [-1,1]$ by the transformation
\begin{equation}
t = 2 \left( \frac{s}{S_{o}} \right)^{2} - 1.
\end{equation}
We then find that
\begin{equation}
\partial_{t} = \frac{S_{o}^{2}}{4} \frac{1}{s} \, \partial_{s}, \hspace{4ex} \dd t = \frac{4}{S_{o}^{2}} s \dd s
\end{equation}
and the transformed volume element reads
\begin{equation}
\dd V = \frac{S_{o}^{2}}{4} h(t) \dd t \dd \phi \dd \eta.
\end{equation}
To proceed, we factor our height function into the following:
\begin{equation}
h(t) = (1-t)^{\frac{1}{2} \chi_{o}} \htilde(t)^{\chi_{h}}
\end{equation}
where $\htilde$, the \emph{non-vanishing part} of the height function $h$, is necessarily strictly positive in $[-1,1]$.
The power $\chi_{h} \in \left\{ 1, \frac{1}{2} \right\}$ denotes whether the polynomial $\htilde$ appears in $h$ as $\sqrt{\htilde}$, corresponding
to $\chi_{h} = \frac{1}{2}$, or simply $\htilde$, corresponding to $\chi_{h} = 1$.
The power $\chi_{o}$ denotes an \emph{outer equatorial singularity} in the domain.
It equals unity when there is a $\sqrt{1-t}$ factor in the height function $h(t)$, otherwise it is zero.

\begin{figure}[H]
\centering
\begin{tabular}{ccc}
    Paraboloid & Oblate Spheroid & Biconcave Disk \\
    \includegraphics[width=0.32\linewidth]{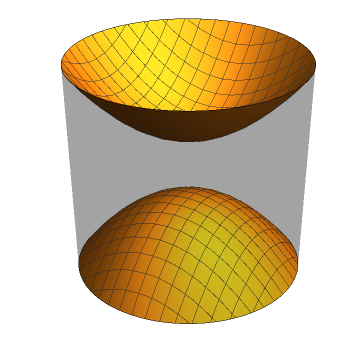} &
    \includegraphics[width=0.32\linewidth]{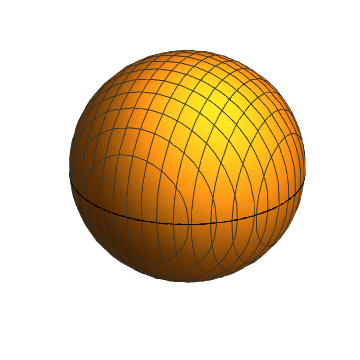}   &
    \includegraphics[width=0.32\linewidth]{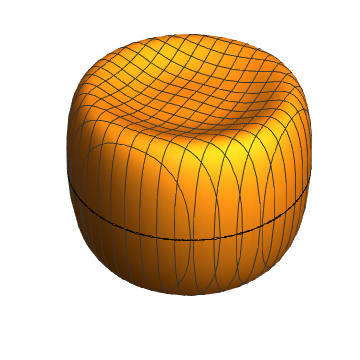} \\
    \parbox{3cm}{\begin{align*}
        \htilde(t) &= \frac{1}{4}(2+t) \\
        h(t)       &= \htilde(t) \\
        \left( \chi_{o}, \chi_{h} \right) &= \left(0, 1\right)
    \end{align*}}
    &
    \parbox{3cm}{\begin{align*}
        \htilde(t) &= 0.8 \\
        h(t)       &= \htilde(t) \, \sqrt{1-t} \\
        \left( \chi_{o}, \chi_{h} \right) &= \left(1, 1\right)
    \end{align*}}
    &
    \parbox{3cm}{\begin{align*}
        \htilde(t) &= \frac{1}{9} \left( 2(1+t)^{2}+1 \right) \\
        h(t)       &= \sqrt{ \htilde(t) } \, \sqrt{1-t} \\
        \left( \chi_{o}, \chi_{h} \right) &= \left(1, \frac{1}{2}\right)
    \end{align*}}
\end{tabular}
\caption{Stretched cylindrical domains.  Paraboloid (left), oblate spheroid (middle) and biconcave disk (right).
Each of these geometries utilizes a distinct height function and combinations of parameters.  The first non-trivial
extension to a standard right cylinder is the paraboloid domain: the height function is a linear, non-vanishing function of $t$.
Note this is a quadratic factor in $s$; geometries linear in $s$ have a cusp at $s = 0$, making the domain non-smooth.
The next extension in weight function is the oblate spheroid, which has a constant $\htilde$ but an outer equatorial singularity at $t = 1$.  Boundary singularities are
always of Jacobi type - they occurs at the edge of the numerical domain and so map onto the $\left( 1 \pm t \right)$ polynomial in the Jacobi weight function.
In the case of the outer equatorial singularity the Jacobi exponent is $a = \frac{1}{2}$.  Cylinders cannot have vanishing height at $s = 0$ and so cannot have
an inner equatorial singularity.
The biconcave disk extends this geometry to multiply by the square root of the non-vanishing part of the height $\htilde(t)$.
This takes the Jacobi $c$ parameter from $1$ to $\frac{1}{2}$.  The gyroscopic polynomial method applies to any of these types of geometries and discretizes
the domain using particular products of vertical and radial polynomials.
}
\label{fig:cylinders}
\end{figure}

Figure~\ref{fig:cylinders} demonstrates three different sample geometries in the stretched
cylindrical coordinate system.  In this figure, the paraboloid has $(\chi_{o},\chi_{h}) = (0,1)$ since there is no outer equatorial singularity ($\chi_{o} = 0$) and the
height function is not square rooted ($\chi_{h} = 1$).  The oblate spheroid has $(\chi_{o},\chi_{h}) = (1,1)$ due to the outer equatorial singularity.
The biconcave disk has $(\chi_{o},\chi_{h}) = \left(1, \frac{1}{2}\right)$ due to both the outer equatorial singularity and the square root non-vanishing height.

\subsection{Stretched Annular Coordinates}\label{sec:stretched_annular_coordinates}
\begin{figure}[H]
\centering
\begin{tabular}{ccc}
    Paraboloid & Sphere & Torus \\
    \includegraphics[width=0.32\linewidth]{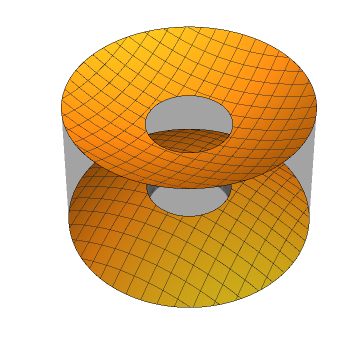} &
    \includegraphics[width=0.32\linewidth]{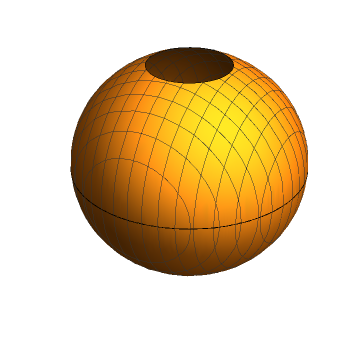}   &
    \includegraphics[width=0.32\linewidth]{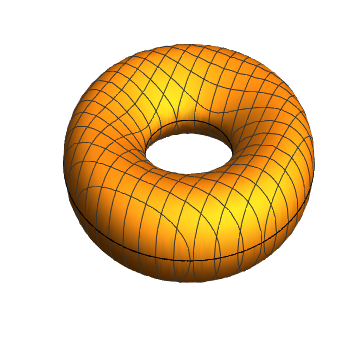} \\
    \parbox{3cm}{\begin{align*}
        h(s)       &= \frac{1}{5}\left(1+s^{2}\right) \\
        \left( \chi_{o}, \chi_{i}, \chi_{h} \right) &= \left(0, 0, 1\right)
    \end{align*}}
    &
    \parbox{3cm}{\begin{align*}
        h(s)       &= \sqrt{1-s^{2}} \\
        \left( \chi_{o}, \chi_{i}, \chi_{h} \right) &= \left(1, 0, 1\right)
    \end{align*}}
    &
    \parbox{3cm}{\begin{align*}
        h(t)       &= \sqrt{1-t^{2}} \\ 
        \left( \chi_{o}, \chi_{i}, \chi_{h} \right) &= \left(1, 1, 0\right)
    \end{align*}}
\end{tabular}
\caption{Stretched annular domains for $s \in \left[S_{i}, S_{o}\right] = \left[0.35, 1\right]$.  Paraboloid (left), sphere (middle). and torus (right).
Similar to Figure~\ref{fig:cylinders}, each of these geometries utilizes a distinct height function and combinations of parameters.
The outer equatorial singularity at $s = S_{o}$ ($t = +1$) leads to a $\sqrt{1-t}$ factor in the height function, represented by $\chi_{o} = 1$,
while the inner $s = S_{i}$ ($t = -1$) equatorial singularity contributes a $\sqrt{1+t}$ factor, corresponding to $\chi_{i} = 1$.  Though we don't
display an example in this figure, gyroscopic polynomials admit square roots around the non-vanishing height $\htilde$ denoted by $\chi_{h} = \frac{1}{2}$.
}
\label{fig:annuli}
\end{figure}

The \emph{stretched annulus coordinate system} ${ (s, \phi, \eta) }$ matches the stretched cylinder except the limited radial domain $s \in [S_{i}, S_{o}]$ - see Figure~\ref{fig:annuli}.
We define the variable $t \in [-1,1]$ by
\begin{equation}
t = \frac{1}{S_{o}^{2} - S_{i}^{2}} \left( 2 s^{2} - \left( S_{o}^{2} + S_{i}^{2} \right) \right).
\end{equation}
Derivatives transform as
\begin{equation}
\partial_{t} = \frac{S_{o}^{2} - S_{i}^{2}}{4} \frac{1}{s} \, \partial_{s}, \hspace{4ex} \dd t = \frac{4}{S_{o}^{2} - S_{i}^{2}} s \dd s
\end{equation}
and the volume element reads
\begin{equation}
\dd V = \frac{S_{o}^{2} - S_{i}^{2}}{4} h(t) \dd t \dd \phi \dd \eta.
\end{equation}
In both the cylindrical and annular cases we changed variables from $s$ to $t \in [-1,1]$.  This transformation replaced the
$s \dd s$ part of the geometric volume element with $\text{d} t$.  This has two benefits: (i) the volume element is simplified
to remove the non-Jacobi-type weight factor $s$ and (ii) the choice of variables guarantees regularity at coordinate singularities.
We describe regularity with respect to coordinate singularities in more detail in Section~\ref{sec:basis}.

Following the notation of Section~\ref{sec:stretched_cylindrical_coordinates} we again factor the height function, this time with a possible zero at the inner radius $t = -1$:
\begin{equation}
h(t) = (1-t)^{\frac{1}{2} \chi_{o}} \, (1+t)^{\frac{1}{2} \chi_{i}} \htilde(t)^{\chi_{h}},
\end{equation}
where $\chi_{o}$ and $\chi_{i}$ denote the presence of the \emph{outer} and \emph{inner equatorial singularities}, respectively, while
$\chi_{h}$ permits the presence of a square root around $\htilde$.
Figure ~\ref{fig:annuli} demonstrates a few particular geometries the coordinate system can describe.

\subsection{Upper-Half Geometries}\label{sec:upper_half_geometries}
\begin{figure}[H]
\centering
\begin{tabular}{ccc}
    Paraboloid & Chebyshev Top \\
    \includegraphics[width=0.32\linewidth]{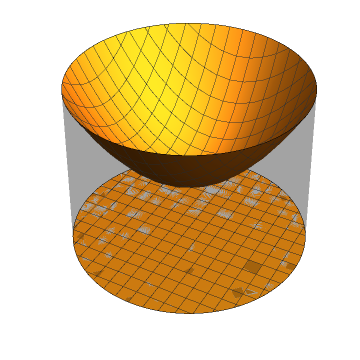} &
    \includegraphics[width=0.32\linewidth]{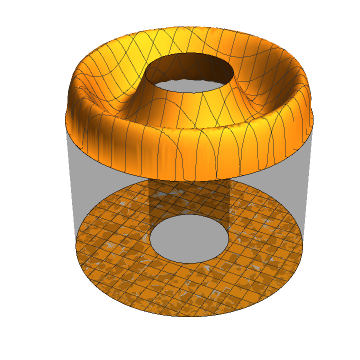} \\
    \parbox{3cm}{\begin{align*}
        h(s) = \htilde(s) &= \frac{1}{8}\left(2+3 s^{2}\right) \\
        \left( \chi_{o}, \chi_{h} \right) &= \left(0, 1\right)
    \end{align*}}
    &
    \parbox{3cm}{\begin{align*}
        h(t) = \htilde(t) &= \frac{1}{16} \left( 8 + 3 t - 4t^{3} \right) \\
        \left( \chi_{o}, \chi_{i}, \chi_{h} \right) &= \left(0, 0, 1\right)
    \end{align*}}
\end{tabular}
\caption{Stretched domains for possible height functions in the upper-half geometry.  This domain only supports $\left(\chi_{o}, \chi_{i}, \chi_{h}\right) = (0,0,1)$.
Note the difference in radial coordinate $s$ (left) and $t$ (right) used to describe the two figures.  We display $h$ expressed in the coordinate that admits the simplest representation.
The right-hand figure displays a more complex height function.  It is proportional to a shifted Chebyshev $T_{3}$ polynomial in the $t$ coordinate and hence is degree six in $s$.
This domain my not be of particular geophysical significance but demonstrates the generality of the gyroscopic method.
}
\label{fig:upperhalf}
\end{figure}

It is often the case that experimenters pursue geometries with flat bottoms, for example tanks with a free surface of
rotating fluid \cite{Lonner_Aggarwal_Aurnou_2022, Cabanes_Aurnou_Favier_Le_Bars_2017,Sommeria_Meyers_Swinney_1989,Lemasquerier_Favier_Le_Bars_2021}
or upper hemispherical rotating flows \cite{Aurnou_Andreadis_Zhu_Olson_2003}.
We have two options for modeling such a domain.
One method is to impose a boundary condition at $z = 0$ ($\eta = 0$) rather than at the numerical domain boundary $z = -h(s)$ ($\eta = -1$).  This boundary
condition at the mid-plane works for geometries with any combination of $(\chi_{o},\chi_{i}) \in \{ 0, 1 \}$ and $\chi_{h} \in \left\{ 1, \frac{1}{2} \right\}$.
In the case that the cylinder or annulus has no equatorial singularities and $\htilde$ appears without a square root - that is, ${ (\chi_{o}, \chi_{i}, \chi_{h}) = (0, 0, 1) }$ - we can shift the domain using a change of coordinates.
We restrict $\eta$ to $[0,1]$ and define $\zeta \in [-1,1]$ by
\begin{equation}
\zeta = 2 \eta - 1.
\iff
\eta = \frac{1}{2} \left( \zeta + 1 \right) 
\end{equation}
This means
\begin{equation}
z = \frac{1}{2} \left( \zeta + 1 \right) h(s).
\end{equation}
This changes not only the scale factor for $z$ derivatives,
\begin{equation}
\dZ = \frac{1}{h(s)} \deta = \frac{2}{h(s)} \dzeta,
\end{equation}
but also the mode coupling in $s$ derivatives,
\begin{equation}
\dS = \ds - \frac{h'(s)}{h(s)} \eta \, \deta = \ds - \frac{h'(s)}{h(s)} \left( 1 + \zeta \right) \, \dzeta.
\end{equation}
As we will see, the structural change $\eta \, \deta \mapsto \left( 1 + \zeta \right) \, \dzeta$ forces us to restrict the upper-half
geometries to non-vanishing polynomials in $t$, namely $h \equiv \htilde$.
Figure~\ref{fig:upperhalf} portrays sample domains in the upper-half geometry.

\subsection{Domain Selection for Exposition}\label{sec:domain_simplification}
For the remainder of the paper we will assume the most basic form of geometry for the cylinder and annulus domains.  This allows us to present the relevant
details on constructing the gyroscopic bases without clouding the presentation with notation and extra indices.  To that end we make the geometric simplifications
$(\chi_{o},\chi_{i},\chi_{h}) = (0,0,1)$ for the remainder of the paper.  This means we have
$h(t) \equiv \htilde(t)$
so that $h$ is a non-vanishing polynomial.  We present the necessary modifications to the gyroscopic basis functions and operators for the special case geometries
in Appendix~\ref{app:general_domains}.

\subsection{Spinor Basis}\label{subsec:spinor_basis}
Though we have fully specified our gyroscopic coordinate systems we still have freedom in how we choose the vector basis to represent vector fields.
The polar geometry of each axial slice leads us to expand vector fields in terms of the \emph{spinor basis}.  One key feature of this basis is that it is
diagonal under differentiation: cross-terms in vector derivatives vanish.  Another important aspect of the spinor basis is that it separates out regularity at the
origin - each component has a different decay rate as $s \to 0$.  By designing these decay rates into the gyroscopic bases we will see the gradient action
on a scalar field maps precisely onto these spin components.  These features of the spinor basis lead to the straightforward generalization to higher rank tensor
fields.  This means we can automatically produce arbitrary rank gradients free of cross-terms; this would be messy in the standard cylindrical vector basis $(\vec{\hat{e}}_{S}, \vec{\hat{e}}_{\Phi}, \vec{\hat{e}}_{Z})$, to say the least. 

We denote the \emph{spin weight} of a vector field component by ${ \sigma \in \{ -1, 0, +1 \} }$.
We define the \emph{spinor basis} following \cite{Vasil_Burns_Lecoanet_Olver_Brown_Oishi_2016,Vasil_Lecoanet_Burns_Oishi_Brown_2019,Ellison_Julien_Vasil_2022}:
\begin{equation}
\vec{\hat{e}}_{\pm} \triangleq \frac{1}{\sqrt{2}} \left(\vec{\hat{e}}_{S} \mp i \vec{\hat{e}}_{\Phi} \right),
\hspace{4ex} \vec{\hat{e}}_{0} \triangleq \vec{\hat{e}}_{Z}
\label{eqn:spinor_basis}
\end{equation}
We represent a vector field $\vec{u}(s,\phi,\eta)$ as
\begin{equation}
\vec{u}(s,\phi,\eta) = u_{S} \vec{\hat{e}}_{S} + u_{\Phi} \vec{\hat{e}}_{\Phi} + u_{Z} \vec{\hat{e}}_{Z} = \sum_{\sigma} u_{\sigma} \vec{\hat{e}}_{\sigma},
\end{equation}
where $u_{S}$, $u_{\Phi}$, $u_{Z}$ and $u_{\sigma}$ are all functions of the stretched coordinates $(s, \phi, \eta)$.
This basis diagonalizes the gradient connection:
\begin{equation}
\nabla_{\pm} \vec{\hat{e}}_{+} = \mp \frac{1}{S} \vec{\hat{e}}_{+}, \hspace{4ex}
\nabla_{\pm} \vec{\hat{e}}_{-} = \pm \frac{1}{S} \vec{\hat{e}}_{-},
\end{equation}
where $\nabla_{\pm}$ acts on a single azimuthal mode ${ e^{i m \phi} f_m(s, \eta) }$ by
\begin{equation}
\nabla_{\pm} \equiv \frac{1}{\sqrt{2}} \left( \partial_{S} \mp \frac{m}{S} \right).
\end{equation}
This diagonalization improves sparsity of differential operators when using the gyroscopic bases to discretize PDEs.

Another important motivation for the spinor basis is that the spin components behave predictably at the $z$ axis.
For ${ \sigma \in \{ -1, 0, +1 \} }$, the $u_{\sigma}$ component of the vector field decays like
\begin{equation}
u_{\sigma}  = \vec{\hat{e}}_{\sigma}^{\dagger} \cdot \vec{u} = \vec{\hat{e}}_{\sigma}^{*} \cdot \vec{u} \sim s^{\left| m \right| + \sigma } F_{\sigma} \left( s^{2} \right), \hspace{4ex} s \to 0,
\end{equation}
where $F_{\sigma}$ is an arbitrary well-behaved function of $s^{2}$.  This means the the basis functions $\Phi_{m,l,k}^{(\alpha, \sigma)}$ previously defined
behave precisely as needed to represent vector fields regular throughout the domain.

Each additional rank $r$ in a tensor field contributes an additional value in $\sigma$'s range.  For example, a rank two tensor $T$ will have a component
\begin{equation}
\vec{\hat{e}}_{+}^{\dagger} \cdot \vec{\hat{e}}_{+}^{\dagger} \cdot T \propto s^{|m| + 2}, \hspace{4ex} s \to 0
\end{equation}
so we see $|\sigma| \le r$.  For more details on the spinor basis for higher rank tensors see \cite{Vasil_Burns_Lecoanet_Olver_Brown_Oishi_2016}.

\section{Generalized Jacobi Polynomials}\label{sec:generalized_jacobi}
We defined the stretched cylindrical and annular coordinate systems to create a numerical domain for computation of differential operators advantageous to gyroscopic flows.
In both cases we mapped the curved upper surface $z = h(s)$ of the domain onto a flat one $\eta = 1$.  To proceed we create a Cartesian $(x,y,z)$ polynomial basis
expressed in the $(t,\phi,\eta)$ coordinates.  The volume element in the stretched coordinates leads to a standard Jacobi weight function in the inner product for
the vertical polynomials but a generalized one for the radial polynomials.  This generalized weight consists of additional non-vanishing polynomial factors.

In this section we develop orthonormal polynomials and the corresponding polynomial algebra with respect to a generalized Jacobi weight function.
Define the \emph{generalized Jacobi weight} with $n$ additional polynomial factors by
\begin{equation}
w^{\left( a,b,c_1,\hdots,c_n; \, p_1,\hdots,p_n \right)}(z) \triangleq (1-z)^{a} (1+z)^{b} p_{1}(z)^{c_{1}} \hdots p_{n}(z)^{c_{n}},
\end{equation}
where $a,b > -1$, $c_{i} \in \mathbb{R}$ and each $p_{i}(z)$ is a non-vanishing polynomial for $z \in \left[-1,1\right]$.  We refer to the polynomials $p_{i}(z)$ as the \emph{augmenting factors} of the
generalized Jacobi weight function.  To ease notation we define the augmenting parameter vector $\vec{c} = \left( c_{1}, \hdots, c_{n} \right)$ 
and augmenting polynomial vector $\vec{p} = \left( p_{1}, \hdots, p_{n} \right)$.

We define the inner product
\begin{equation}
\left \langle f, g \right \rangle_{\left(a,b,\vec{c};\,\vec{p}\right)}
    \triangleq \int_{-1}^{1} \dd z \, w^{\left(a,b,\vec{c}; \, \vec{p}\right)}(z) \, \overline{f} g
    = \int_{-1}^{1} \dd z \, (1-z)^{a} (1+z)^{b} p_{1}(z)^{c_{1}} \cdots p_{n}(z)^{c_{n}} \, \overline{f} g
\label{eqn:generalized_jacobi_inner_product}
\end{equation}
and corresponding norm
\begin{equation}
\norm{ f }_{\left(a,b,\vec{c};\,\vec{p}\right)}^{2} \triangleq \left\langle f, f \right\rangle_{\left(a,b,\vec{c};\,\vec{p}\right)}.
\end{equation}
The inner product (\ref{eqn:generalized_jacobi_inner_product}) induces family of Hilbert spaces $\mathcal{H}\left(a,b,\vec{c};\,\vec{p}\right)$, where
\begin{equation}
\mathcal{H}\left(a,b,\vec{c};\,\vec{p}\right) \triangleq \left\{ f : [-1,1] \to \mathbb{C} : \norm{ f }_{\left(a,b,\vec{c};\,\vec{p}\right)} < \infty \right\}.
\end{equation}
The generalized Jacobi polynomials $P_{n}^{\left(a,b,\vec{c};\,\vec{p}\right)}$ are those polynomials orthonormal under 
the weight function $w^{\left( a,b,\vec{c};\,\vec{p} \right)}$:
\begin{equation}
\left \langle P_{n}^{\left( a,b,\vec{c};\,\vec{p} \right)}, P_{m}^{\left( a,b,\vec{c};\,\vec{p} \right)} \right \rangle_{\left(a,b,\vec{c};\,\vec{p}\right)}
    = \delta_{m,n}
\end{equation}
Normalization requires the mass of the inner product:
\begin{equation}
\omega^{\left(a,b,\vec{c};\,\vec{p}\right)} \triangleq \norm{ 1 }_{\left(a,b,\vec{c};\,\vec{p}\right)}
\end{equation}
so that
\begin{equation}
P_{0}^{\left(a,b,\vec{c};\,\vec{p}\right)}(z)  = \frac{1}{\sqrt{\omega^{\left(a,b,\vec{c};\,\vec{p}\right)}}}.
\end{equation}
From here on we drop the $\vec{p}$ symbol from our notation when it is unambiguous.  Operations on the polynomials never change the underlying
augmenting polynomial factors, just their powers, so for notational simplicity we retain only the indices $\left(a,b,\vec{c}\right)$.  When our
discussion doesn't depend on Jacobi the parameters $a,b$ we also drop these.
Computing these polynomials requires some additional tools \cite{Gautschi_1982,Magnus_1993,Olver_Townsend_Vasil_2018,Olver_Townsend_Vasil_2019,Olver_Slevinsky_Townsend_2020,Snowball_Olver_2020a,Snowball_Olver_2020b}, but they all boil down to using
known quadrature rules to compute the three-term recurrence coefficients $\left( \alpha_{n}^{\left(a,b,\vec{c}\right)}, \beta_{n}^{\left(a,b,\vec{c}\right)} \right)$ satisfying
\begin{equation}
z P_{n}^{\left(a,b,\vec{c}\right)}(z) = \beta_{n}^{\left(a,b,\vec{c}\right)} P_{n+1}^{\left(a,b,\vec{c}\right)}(z) + \alpha_{n}^{\left(a,b,\vec{c}\right)} P_{n}^{\left(a,b,\vec{c}\right)}(z) + \beta_{n-1}^{\left(a,b,\vec{c}\right)} P_{n-1}^{\left(a,b,\vec{c}\right)}(z).
\label{eqn:three_term_recurrence}
\end{equation}

We describe a method for computing $(\alpha_{n}, \beta_{n})$ in Section~\ref{sec:jacobi_three_term_recurrence}.
Once we have the three-term recurrence coefficients in hand we can compute the polynomials themselves.
Importantly, the recurrence (\ref{eqn:three_term_recurrence}) also enables us to compute exact quadrature rules for polynomials up to a specified degree.
We describe the calculation in Section~\ref{sec:jacobi_quadrature_rule}.
This quadrature is the final tool necessary to compute the sparse operator matrix coefficients acting on the polynomial bases.
We define embedding operators in Section~\ref{sec:jacobi_embedding_operators} and differential operators in Section~\ref{sec:jacobi_differential_operators},
then describe how to compute them in Section~\ref{sec:jacobi_operator_computation}.

\subsection{Three-term Recurrence}\label{sec:jacobi_three_term_recurrence}
The three-term recurrence (\ref{eqn:three_term_recurrence}) is critical to computing with orthogonal polynomial systems.  From these coefficients we can compute
the polynomials and their derivatives as well as quadrature rules and the matrix entries for discretized representations of operators acting on the polynomial systems.  We need all
of these tools before we can build our gyroscopic polynomial bases.

There are many methods for computing the three-term recurrence
coefficients for orthogonal polynomials with respect to a given measure \cite{Wheeler_1974,Gautschi_1982}, including the discretized Stieltjes procedure
and modified Chebyshev algorithm.  These methods use modified moments of the measure computed with known Gaussian quadrature rules
to calculate the recurrence coefficients.
As we will see, vertical degree parameterizes the radial parameters of the gyroscopic polynomials.  This means one or more of the generalized Jacobi parameters
increments with this degree.  To generate our recurrence coefficients we therefore turn to the Christoffel-Darboux formulation \cite{Snowball_Olver_2020a}
to recursively generate three-term coefficients for orthogonal polynomials with successively incremented parameter values.
Below, we adapt the formulation \cite{Snowball_Olver_2020a} to the interval $[-1,1]$ with a slightly more general augmenting polynomial.
The proof in \cite{Snowball_Olver_2020a} using the Christoffel-Darboux formula only needs minor modification; we omit it for brevity.

Let $w^{(0)}(z)$ be a weight function with known three-term recurrence coefficients $\left(\alpha_{n}^{(0)}, \beta_{n}^{(0)}\right)$ and orthonormal polynomials $P_{n}^{(0)}(z)$.
In practice this system is a Jacobi polynomial system with fixed $a,b$ parameters, but since these parameters are fixed throughout the discussion we omit them here.  The goal is to 
derive recurrence coefficients corresponding to augmenting the weight function with a linear polynomial.
To that end, for positive integer $c$ define $w^{(c)}(z) = w^{(0)}(z) \, p(z)^{c}$ where $p(z) = m \, z + b$ is a linear polynomial in $z$.  We want to compute $\left(\alpha^{(c)}_{n}, \beta^{(c)}_{n}\right)$, the
three-term recurrence coefficients for $w^{(c)}(z)$, in terms of known quantities.

Let $z_{0} = -\frac{b}{m}$ be the $z$ intercept of $p(z)$.  For integer $0 \le k < c$ define $C_{n}^{(k)}$ such that
\begin{equation}
C_{n}^{(k)} = \left(\pm 1 \right)^{n} \sqrt{ \frac{\pm 1}{ P_{n}^{(k)}\left(z_{0}\right) P_{n+1}^{(k)}\left(z_{0}\right) \beta_{n}^{(k)} } },
\end{equation}
where we select $+1$ when $m < 0$ and $-1$ otherwise.  Then
\begin{equation}
\begin{aligned}
\alpha_{n}^{(k+1)} &= \frac{P_{n+2}^{(k)}\left(z_{0}\right)}{P_{n+1}^{(k)}\left(z_{0}\right)} \beta_{n+1}^{(k)} - \frac{P_{n+1}^{(k)}\left(z_{0}\right)}{P_{n}^{(k)}\left(z_{0}\right)} \beta_{n}^{(k)} + \alpha_{n+1}^{(k)}, \\
\beta_{n}^{(k+1)}  &= \frac{C_{n}^{(k)}}{C_{n+1}^{(k)}} \frac{P_{n}^{(k)}\left(z_{0}\right)}{P_{n+1}^{(k)}\left(z_{0}\right)} \beta_{n}^{(k)}
\end{aligned}
\end{equation}
and
\begin{equation}
P_{n}^{(k+1)}\left(z_{0}\right) = C_{n}^{(k)} \sum_{j=0}^{n} \left( P_{j}^{(k)}\left(z_{0}\right) \right)^{2}.
\end{equation}
This formulation provides a means to recursively compute the coefficients for $w(z) = w_{0}(z) p(z)^{c}$ in terms of the weight function $w_{0}(z) p(z)^{c-1}$.
We require an initialization of $\left(\alpha_{n}^{(0)}, \beta_{n}^{(0)}\right)$ of size $N + c + 1$ to compute $N$ terms of $\left(\alpha_{n}^{(c)}, \beta_{n}^{(c)}\right)$.

The Christoffel-Darboux method applies to linear factors in the weight function.  For weight functions with higher
degree polynomials we factor the polynomials into a product of linear and quadratic factors.  A bit of algebra allows us to bypass
complex-valued weight functions for complex-conjugate pairs of roots.  We show the generalization to quadratic factors with complex conjugate
roots in Appendix~\ref{app:christoffel_darboux}.  The Christoffel-Darboux formulation therefore applies to any real-coefficient augmenting polynomial.

To compute the three-term recurrence for the generalized Jacobi weight $w^{\left(a,b,\vec{c};\,\vec{p}\right)}(z)$ with $n$ augmenting factors we repeat this process for each factor $p_{i}$ individually.
At each stage we replace $w^{(0)}(z)$ with the most recently computed polynomial system.
This increases the corresponding initialization size of the known three-term recurrence coefficients but otherwise the process in unchanged.

In what follows we see the gyroscopic bases are built in a hierarchy of generalized Jacobi parameters.  Computing the three-term recurrence for each member of the hierarchy from scratch
using the Stieltjes or modified Chebyshev procedures requires $\mathcal{O}\left(N^{2}\right)$ operations, or $\mathcal{O}\left(N^{3}\right)$ total operations for the entire hierarchy. 
The Christoffel-Darboux formulation enables computing the same set of three-term recurrences in $\mathcal{O}\left(N^{2}\right)$ total operations, making the gyroscopic
basis computationally tractable.

\subsection{Quadrature Rules}\label{sec:jacobi_quadrature_rule}
Starting with the three-term recurrence $(\alpha_{n}, \beta_{n})$ for a general polynomial system with $P_{n}$ orthonormal with respect to weight function $w(z)$ we compute the $N$-term Gauss quadrature rule $(z_{j}, w_{j})$ such that
\begin{equation}
\int_{-1}^{1} \dd z \, w(z) f(z) \approx \sum_{j=0}^{N-1} w_{j} \, f(z_{j}),
\end{equation}
where we have equality when $f$ is a polynomial of degree at most $2N-1$.  The Golub-Welsch algorithm \cite{Olver_Slevinsky_Townsend_2020} for computing $(z_{j}, w_{j})$ utilizes the operator
$\mathcal{Z} P_{n} = z P_{n}(z)$, which forms a symmetric tridiagonal matrix $Z$ acting on the basis $P_{n}$ with diagonal $\alpha_{n}$ and off-diagonals $\beta_{n}$.
We truncate (\ref{eqn:three_term_recurrence}) to degree $N$ and rewrite it in matrix form via
\begin{equation}
Z_{0:N-1,0:N} \begin{pmatrix} P_{0}(\lambda) \\ \vdots \\ P_{N}(\lambda) \end{pmatrix} = \lambda \begin{pmatrix} P_{0}(\lambda) \\ \vdots \\ P_{N-1}(\lambda) \end{pmatrix}
\end{equation}
so that
\begin{equation}
\left( Z_{N} - \lambda I \right) \begin{pmatrix} P_{0}(\lambda) \\ \vdots \\ P_{N-2}(\lambda) \\ P_{N-1}(\lambda) \end{pmatrix} = \begin{pmatrix} 0 \\ \vdots \\ 0 \\ -\beta_{N-1} P_{N}(\lambda) \end{pmatrix}.
\end{equation}
Now the right-hand side is the zero vector only when $\lambda$ is a root of $P_{N}(\lambda)$.  The eigenvalues $z_j$ are therefore the quadrature nodes,
and the associated eigenvector is then 
\begin{equation}
\begin{pmatrix} P_{0}(z_j) \\ \vdots \\ P_{N-1}(z_j) \end{pmatrix}
\end{equation}
We compute the quadrature weights via
\begin{equation}
\frac{1}{w_{j}} = \sum_{i=0}^{N-1} P_{i}(z_{j})^{2}.
\end{equation}
The algorithm takes $\mathcal{O}(N^{2})$ operations to compute due to the symmetric tridiagonal structure of the matrix $Z$.
In case the eigensolve has error beyond machine precision we refine the eigenvalue solution to $z_j$ using a Newton iteration,
\begin{equation}
z_{j}^{(k+1)} = z_{j}^{(k)} - \frac{ P_{j}\left(z_{j}^{(k)}\right) }{ P_{j}^{\prime}\left(z_{j}^{(k)}\right) },
\end{equation}
for $0 < k < \infty$.
We then recompute $w_j$ using the final iterate $z_{j}^{(\infty)}$.  In practice we need only a single iteration to converge within machine precision.

\subsection{Embedding Operators}\label{sec:jacobi_embedding_operators}
Just like for classical Jacobi polynomials, the generalized Jacobi polynomials have sparse embedding operators.  We define
\begin{equation}
\begin{aligned}
\mathcal{I}_{a} &: \mathcal{H}\left(a,b,\vec{c}\right) \mapsto \mathcal{H}\left(a+1,b,\vec{c}\right) \\
    &\hspace{4ex}\mathcal{H}\left(a,b,\vec{c}\right) \ni f =  \mathcal{I}_{a} f \in \mathcal{H}\left(a+1,b,\vec{c}\right) \\
\mathcal{I}_{b} &: \mathcal{H}\left(a,b,\vec{c}\right) \mapsto \mathcal{H}\left(a,b+1,\vec{c}\right) \\
    &\hspace{4ex}\mathcal{H}\left(a,b,\vec{c}\right) \ni f =  \mathcal{I}_{b} f \in \mathcal{H}\left(a,b+1,\vec{c}\right) \\
\mathcal{I}_{c_i} &: \mathcal{H}\left(a,b,c_1,\hdots,c_i,\hdots,c_n\right) \mapsto \mathcal{H}\left(a,b,c_1,\hdots,c_i+1,\hdots,c_n\right) \\
    &\hspace{4ex}\mathcal{H}\left(a,b,c_1,\hdots,c_i,\hdots,c_n\right) \ni f = \mathcal{I}_{c_{i}} f \in \mathcal{H}\left(a,b,c_1,\hdots,c_i+1,\hdots,c_n\right)
\end{aligned}
\end{equation}
where each operator is the identity acting to embed a function in the corresponding codomain.  This implies the Hilbert spaces naturally nest:
for $\delta_{a},\delta_{b},\delta_{c_{i}} \in \mathbb{N}$ we have
\begin{equation}
\mathcal{H}\left(a,b,\vec{c}\right) \subseteq \mathcal{H}\left(a+\delta_{a},b+\delta_{b},\vec{c}+\vec{\delta_{c}}\right).
\end{equation}

We can achieve multiplication by the various weight factors using the embedding adjoints.  We have
\begin{equation}
\begin{aligned}
\mathcal{I}_{a}^{\dagger} &: \mathcal{H}\left(a,b,\vec{c}\right) \mapsto \mathcal{H}\left(a-1,b,\vec{c}\right) \\ &\hspace{4ex} \mathcal{I}_{a}^{\dagger} f = (1-z) f \\
\mathcal{I}_{b}^{\dagger} &: \mathcal{H}\left(a,b,\vec{c}\right) \mapsto \mathcal{H}\left(a,b-1,\vec{c}\right) \\ &\hspace{4ex} \mathcal{I}_{b}^{\dagger} f = (1+z) f \\
\mathcal{I}_{c_i}^{\dagger} &: \mathcal{H}\left(a,b,c_1,\hdots,c_i,\hdots,c_n\right) \mapsto \mathcal{H}\left(a,b,c_1,\hdots,c_i-1,\hdots,c_n\right) \\ &\hspace{4ex} \mathcal{I}_{c_{i}}^{\dagger} f = p_{i}(z) f.
\end{aligned}
\end{equation}
The bandwidth of each operator is one more than the degree of the corresponding weight factor.

\begin{figure}[H]
\centering
\includegraphics[width=\linewidth]{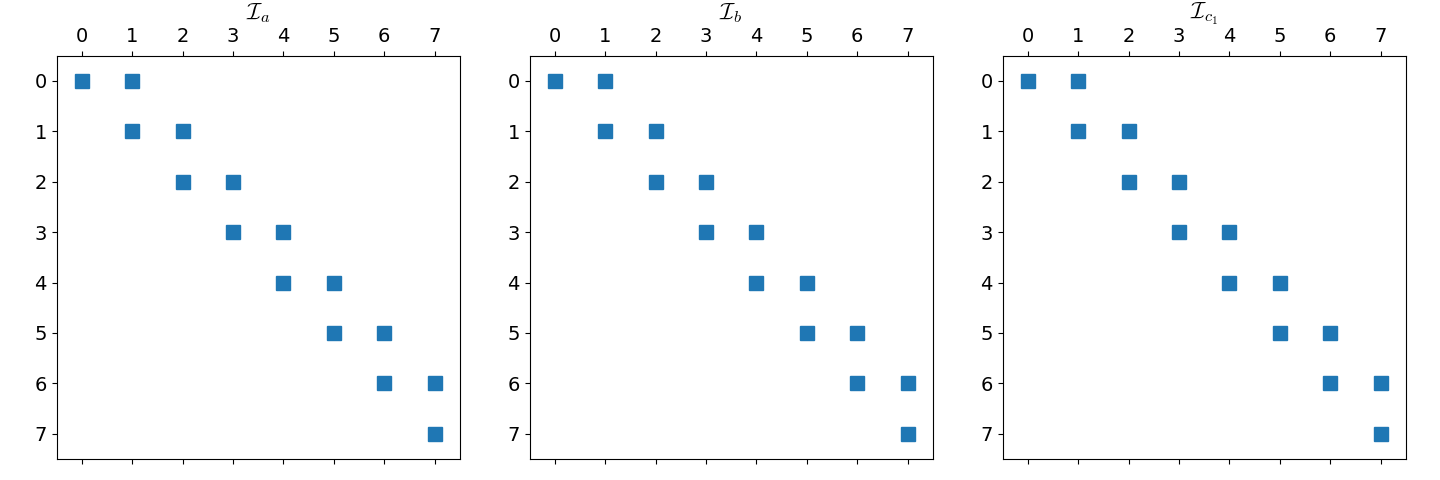}
\caption{Sparsity diagram for generalized Jacobi embedding operators for weight function with a single additional linear factor.
We interpret row $m$ of the operator as recombining the expansion coefficients expressed in domain basis $P$ to form the coefficient for codomain polynomial $Q_{m}$.
We interpret column $n$ of the operator as recombining the codomain polynomials $Q$ to produce $\mathcal{L} P_{n}$.
The adjoints of each operators have a transposed sparsity structure.}
\label{fig:genjacobi_sparsity_embed}
\end{figure}

Let a field be represented by a column vector of expansion coefficients:
\begin{equation}
f(z) = \begin{pmatrix} P_{0}(z) & \hdots & P_{n}(z) \end{pmatrix} \begin{pmatrix} \widehat{F}_{0} \\ \vdots \\ \widehat{F}_{n} \end{pmatrix}.
\end{equation}
Then an operator $\mathcal{L}$ acting on domain with orthogonal polynomial system $\left\{ P_{n} \right\}$ and codomain with orthogonal polynomial system $\left\{ Q_{m} \right\}$ takes the matrix form $L_{m,n}$ such that
\begin{equation}
\mathcal{L} f(z) = \begin{pmatrix} Q_{0}(z) & \hdots & Q_{m}(z) \end{pmatrix} 
    \begin{pmatrix} L_{0,0} & \hdots & L_{0,n} \\ \vdots & \ddots & \vdots \\ L_{m,0} & \hdots & L_{m,n} \end{pmatrix}
    \, \begin{pmatrix} \widehat{F}_{0} \\ \vdots \\ \widehat{F}_{n} \end{pmatrix}.
\end{equation}
We then see that row $m$ of $L$ recombines the expansion coefficients $\vec{\widehat{F}}$ expressed in the $P$ basis to form the coefficient of $\mathcal{L} f$ multiplying $Q_{m}$.
Equivalently, we can interpret $L$ acting on the basis functions themselves:
\begin{equation}
\mathcal{L} \vec{P}(z) = \begin{pmatrix} \mathcal{L} P_{0}(z) & \hdots & \mathcal{L} P_{n}(z) \end{pmatrix} = \begin{pmatrix} Q_{0}(z) & \hdots & Q_{m}(z) \end{pmatrix} 
    \begin{pmatrix} L_{0,0} & \hdots & L_{0,n} \\ \vdots & \ddots & \vdots \\ L_{m,0} & \hdots & L_{m,n} \end{pmatrix},
\end{equation}
in which case column $n$ of $L$ recombines the $Q$ basis to form $\mathcal{L} P_{n}$.  These points of view help us interpret
Figure~\ref{fig:genjacobi_sparsity_embed}, where plot the sparsity structure of the embedding operators for a single first order augmenting polynomial.

\subsection{Differential Operators}\label{sec:jacobi_differential_operators}
A single framework unifies all sparse differential operators on the generalized Jacobi polynomials.  We define the operator
\begin{equation}
\mathcal{D}\left( \delta_{a}, \delta_{b}, \vec{\delta_c} \right)
    : \mathcal{H}\left(a,b,\vec{c}\right) \mapsto \mathcal{H}\left(a+\delta_{a}, b+\delta_{b}, \vec{c}+\vec{\delta_c} \right)
\end{equation}
where each $\delta \in \{-1,+1\}$.  To compute the action of the operator we need some definitions.
Let $k \le n$ and $\left( i_{1}, \hdots, i_{k} \right)$ be a strictly increasing index tuple $1 \le i_{1} < i_{2} < \hdots < i_{k} \le n$.
Define $\rho_{i_{1},\hdots,i_{k}}$ and $\rho_{i_{1},\hdots,i_{k}}'$ such that
\begin{equation}
\begin{aligned}
\rho_{i_{1},\hdots,i_{k}}(z) &\triangleq \prod_{j=1}^{k} p_{i_j}(z) \\
\rho_{i_{1},\hdots,i_{k}}'(z) &\triangleq \sum_{j=1}^{k} c_{i_j} \, p_{i_j}'(z) \underset{l \ne j}{\prod_{l = 1}^{k}} p_{i_l}(z),
\end{aligned}
\end{equation}
where $\rho_{\emptyset}(z) \triangleq 1$ and $\rho_{\emptyset}'(z) \triangleq 0$.
Then our operator takes the following forms:
\begin{equation}
\begin{aligned}
\mathcal{D}\left(+1, +1, \vec{\delta_c} \right) &=
     \rho_{i_{1}, \hdots, i_{k}}(z) \frac{d}{dz} + \rho_{i_{1}, \hdots, i_{k}}'(z) \\
\mathcal{D}\left(+1, -1, \vec{\delta_c}, \right) &=
     \rho_{i_{1}, \hdots, i_{k}}(z) \left(\hphantom{-}b + (1+z) \frac{d}{dz} \right) + \rho_{i_{1}, \hdots, i_{k}}'(z) (1+z) \\
\mathcal{D}\left(-1, +1, \vec{\delta_c}, \right) &=
     \rho_{i_{1}, \hdots, i_{k}}(z) \left(          - a + (1-z) \frac{d}{dz} \right) + \rho_{i_{1}, \hdots, i_{k}}'(z) (1-z) \\
\mathcal{D}\left(-1, -1, \vec{\delta_c}, \right) &=
     \rho_{i_{1}, \hdots, i_{k}}(z) \left(-(1+z)a + (1-z)b + (1-z^2) \frac{d}{dz} \right) + \rho_{i_{1}, \hdots, i_{k}}'(z) (1-z^2)
\end{aligned}
\end{equation}
where the index tuple $(i_j)$ contains index $l$ if $\delta_{c_{l}} = -1$.  We see these operators reduce to the standard Jacobi polynomial operators
when the index tuple is empty.  Hence:
\begin{equation}
\begin{aligned}
\mathcal{D}\left(+1, +1, \vec{+1} \right) &=
     \frac{d}{dz} &&\equiv \mathcal{D}_{m}, \\
\mathcal{D}\left(+1, -1, \vec{+1} \right) &=
     \hphantom{-}b + (1+z) \frac{d}{dz} &&\equiv \mathcal{D}_{s}, \\
\mathcal{D}\left(-1, +1, \vec{+1} \right) &=
     - a + (1-z) \frac{d}{dz} &&\equiv -\mathcal{D}_{s}^{\dagger}, \\
\mathcal{D}\left(-1, -1, \vec{+1} \right) &=
     -(1+z)a + (1-z)b + (1-z^2) \frac{d}{dz} &&\equiv -\mathcal{D}_{m}^{\dagger},
\end{aligned}
\label{eqn:genjacobi_diffops_jacobi}
\end{equation}
with $\mathcal{D}_{m}$ and $\mathcal{D}_{s}$ defined as in Appendix C of \cite{Vasil_Lecoanet_Burns_Oishi_Brown_2019}.

\begin{figure}[H]
\centering
\includegraphics[width=\linewidth]{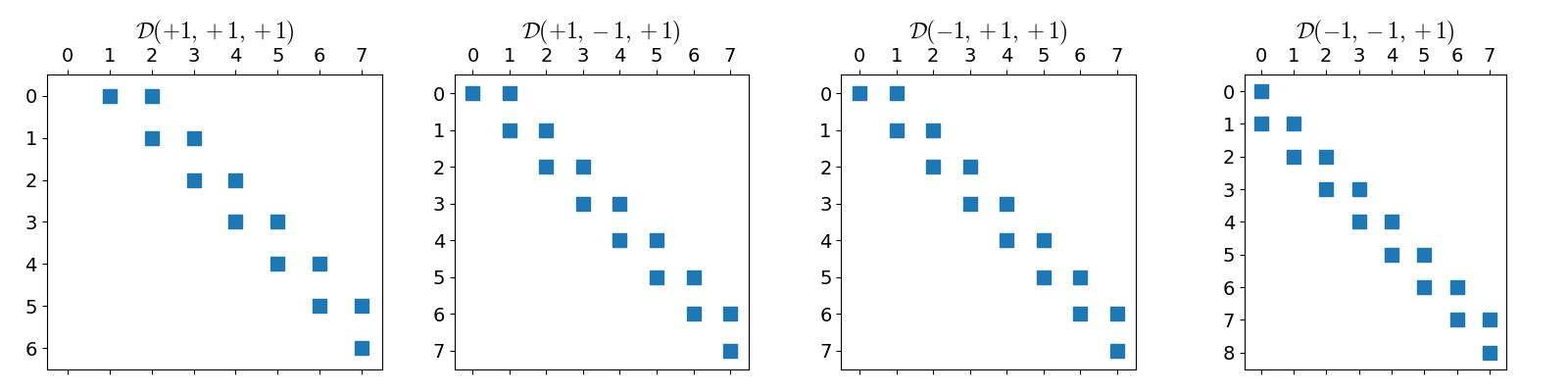}
\caption{Sparsity diagram for generalized Jacobi differential operators (\ref{eqn:genjacobi_diffops_jacobi}) for weight function with a single additional linear factor.
The operators act on column vectors of expansion coefficients from the left, with row $m$ producing the $m^{\text{th}}$ expansion coefficient in the codomain basis.
The adjoint operators have a transposed sparsity structure and swapped domain and codomain compared with the non-adjoint operators.}
\label{fig:genjacobi_sparsity_diffops}
\end{figure}

The bandwidth of the each operator is one more than the degree of the $\rho_{1,\hdots,n}(z)$ polynomial.  We therefore see sparsity depends on the specific $p_{i}(z)$
that appear in the generalized Jacobi weight function.  Each additional degree in the generalized weight contributes an additional diagonal in the operator matrices.
In Figure~\ref{fig:genjacobi_sparsity_diffops} we plot the sparsity structure of the differential operators for a single first order augmenting polynomial.
To provide an example, in the case of stretched cylindrical coordinates with a parabolic height profile we have a single additional factor of degree one.
Differential operators acting on the radial parts of the parabolic cylinder therefore have bandwidth two.  Figure~\ref{fig:genjacobi_sparsity_diffops} shows the sparsity structure for these operators.
The stretched annulus with parabolic height function has two linear augmenting polynomials, and so differential operators have a bandwidth of three.

This sparsity is the key reason we developed the generalized polynomials.  We use the polynomials to build a hierarchy of basis functions orthonormal with respect to the weighted volume element of specific curved geometries.
Incorporating the volume element in the inner product means the basis functions naturally conform to geometric singularities.  Not only are the basis functions well-behaved for representing arbitrary fields in the domain,
but also their nodes don't cluster near coordinate singularities \cite{olver2013fast,Vasil_Burns_Lecoanet_Olver_Brown_Oishi_2016,Lecoanet_Vasil_Burns_Brown_Oishi_2019}.  This has clear benefits for the CFL condition in explicit time-stepping schemes.

\subsection{Operator Computation}\label{sec:jacobi_operator_computation}
The operator matrix entries $L_{m,n}$ for a given operator $\mathcal{L} : \left(a,b,\vec{c}\right) \mapsto \left( a+\delta_{a}, b+\delta_{b}, \vec{c} + \vec{\delta_{c}} \right)$ are
\begin{equation}
L_{m,n} = \left \langle \mathcal{L} P_{n}^{\left(a,b,\vec{c}\right)}, P_{m}^{\left( a+\delta_{a}, b+\delta_{b}, \vec{c} + \vec{\delta_{c}} \right)} \right \rangle_{\left( a+\delta_{a}, b+\delta_{b}, \vec{c}+\vec{\delta_c} \right)}.
\label{eqn:operator_inner_product}
\end{equation}
The three-term recurrence provides a stable method to evaluate the polynomials and their derivatives.  We differentiate (\ref{eqn:three_term_recurrence}) to obtain
\begin{equation}
P_{n+1}'(z) = \frac{1}{\beta_{n}} \left[ (z - \alpha_{n}) P_{n}'(z) + P_{n}(z) - \beta_{n-1} P_{n-1}'(z) \right].
\end{equation}
In Section~\ref{sec:jacobi_quadrature_rule} we computed the quadrature nodes $z_{j}$ and weights $w_{j}$ such that
\begin{equation}
\int_{-1}^{1} \dd z \, w^{\left(a,b,\vec{c}\right)}(z) f(z) \approx \sum_{j=0}^{N-1} w_{j}^{\left(a,b,\vec{c}\right)} f(z_{j}),
\end{equation}
where we have equality for $f$ a polynomial of degree at most $2 N - 1$.
We therefore compute $L_{m,n}$ by evaluating the inner product (\ref{eqn:operator_inner_product}) with the quadrature rule,
using the three-term recurrence to compute the polynomials and their derivatives in the argument of the integral.

\section{The Gyroscopic Bases}\label{sec:basis}
With the tools developed in Section~\ref{sec:generalized_jacobi} we now define a family of basis functions orthogonal in the stretched domain.  Each basis function should be a Cartesian $(x,y,z)$ polynomial expressed in the
stretched coordinate system $(t,\phi,\eta)$.  For this to be the case we note that
\begin{equation}
z^{l} = \eta^{l} h(t)^{l} = \eta^{l} (1-t)^{\frac{l}{2} \chi_{o}} (1+t)^{\frac{l}{2} \chi_{i}} \htilde(t)^{l \, \chi_{h}}
\label{eqn:z_to_the_ell}
\end{equation}
or, using our domain simplification assumption $(\chi_{o},\chi_{i},\chi_{h}) = (0,0,1)$ explained in Section~\ref{sec:domain_simplification},
\begin{equation}
z^{l} = \eta^{l} \htilde(t)^{l}.
\end{equation}
Any gyroscopic polynomial of vertical degree $l$ has an $l$~-dependent prefactor so that it can factor into a Cartesian polynomial.
In what follows we define the gyroscopic basis functions for the stretched cylinder and stretched annulus domains.

\subsection{Stretched Cylinder Bases}
Recall the volume measure $\dd V$ (\ref{eqn:volume_element_cylinder}).  We want to form a three-dimensional basis orthonormal under this measure.
For standard Jacobi polynomials we know we maintain sparsity under differentiation by expressing the output of the derivative with incremented parameters $(a+1,b+1)$.
To that end we weight the volume element $\dd V$ with the boundary polynomial 
\begin{equation}
m^{(\alpha)}(t, \eta) \triangleq \left( 1 - \eta^{2} \right)^{\alpha} \left( 1 - t \right)^{\alpha} \htilde(t)^{2 \alpha}.
\label{eqn:malpha_cylinder}
\end{equation}
This factor multiplying the volume measure gives us access to a numerical parameter $\alpha > -1$ that we let increment under differentiation, a key element to maintaining sparsity of discretized PDEs.
Now define the weighted volume measure $\dd \mu(\alpha)$ by
\begin{equation}
\begin{aligned}
\dd \mu (\alpha) &\triangleq \frac{1}{2 \pi} \frac{4}{S_{o}^{2}} m^{(\alpha)}(t, \eta) \dd V \\
    &= \frac{1}{2 \pi} \left( 1-\eta^{2} \right)^{\alpha} \left( 1 - t \right)^{\alpha} \htilde(t)^{2 \alpha + 1} \dd \phi \dd \eta \dd t
\end{aligned}
\end{equation}
and inner product
\begin{equation}
\left \langle f, g \right \rangle_{\dd \mu(\alpha)} \triangleq \int \dd \mu(\alpha) \overline{f} g.
\end{equation}
This leads us to define the basis function as follows:
\begin{equation}
\Phi_{m,l,k}^{(\alpha,\sigma)} = e^{i m \phi} (1+t)^{\frac{\left|m\right| +\sigma}{2}} \htilde(t)^{l} P_{l}^{(\alpha,\alpha)}(\eta) \, Q_{k}^{\left( \alpha, \left|m\right|+\sigma, 2l + 2\alpha + 1; \htilde \right)}(t).
\label{eqn:cylinder_basis}
\end{equation}
Here $Q_{k}^{(a,b,c;\, p)}$ denotes the degree-$k$ \emph{generalized Jacobi polynomial} orthogonal under the weight $(1-t)^{a} (1+t)^{b} p(t)^{c}$.
These basis functions are orthonormal with respect to $\dd \mu(\alpha)$ since
\begin{equation}
\begin{aligned}
\left \langle \Phi_{m,l,k}^{(\alpha,\sigma)}, \Phi_{m',l',k'}^{(\alpha,\sigma)} \right \rangle_{\dd \mu(\alpha)} 
    &= \frac{1}{2 \pi} \int_{0}^{2 \pi} \dd \phi \, e^{-i \left(m - m'\right) \phi} \int_{-1}^{1} \dd \eta \left(1-\eta^{2}\right)^{\alpha} P_{l}^{(\alpha,\alpha)}(\eta) P_{l'}^{(\alpha,\alpha)}(\eta) \\
    & \times \int_{-1}^{1} \dd t \, \left(1-t\right)^{\alpha} \left( 1 + t \right)^{\left|m\right| + \sigma} \htilde(t)^{ 2 l + 2 \alpha + 1 } \\
    &\hspace{8ex}     Q_{k}^{\left( \alpha, \left|m\right|+\sigma, 2l + 2\alpha + 1; \htilde \right)}(t)
                     Q_{k'}^{\left( \alpha, \left|m\right|+\sigma, 2l + 2\alpha + 1; \htilde \right)}(t) \\
    &= \delta_{m,m'} \, \delta_{l,l'} \, \delta_{k,k'}.
\end{aligned}
\end{equation}
The parameter $\sigma$ is the \emph{spin weight} of the basis as described in Section~\ref{subsec:spinor_basis}.
A key feature of these basis functions is that they conform to any coordinate singularities in the domain.  The full cylinder has a coordinate singularity at $s = 0$.  The basis resolves
this singularity by using the $t$ coordinate, a function of $s^{2}$, in conjunction with the $\left(1+t\right)^{\frac{m}{2}} \propto s^{m}$ prefactor for azimuthal mode $e^{i m \phi}$.
The $\alpha$ parameter creates a mechanism to maintain sparsity under differential operators.  Just as the generalized Jacobi operators increment and decrement their parameters - see Section~\ref{sec:jacobi_differential_operators} -
so will calculus operators applied to the gyroscopic basis functions.  The parameter $\alpha$ enables proper incrementing and decrementing so that the hierarchy of bases fits together in a sparse way.
We must also ensure all terms in an equation live in the same Hilbert space, which we achieve by judicious conversion of the $\alpha$ parameter for each term.

\begin{figure}[H]
\centering
\includegraphics[width=\linewidth]{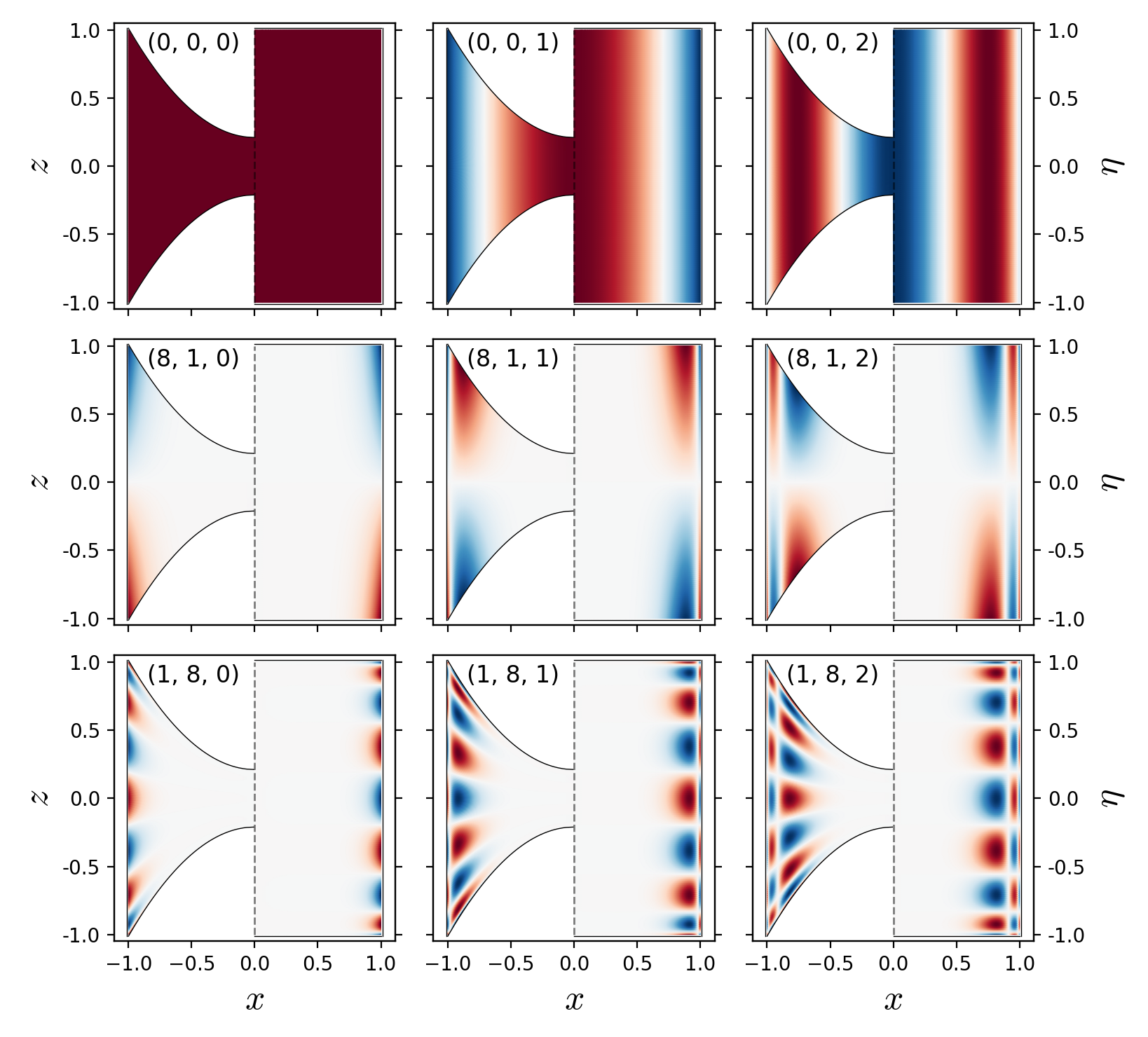}
\caption{Meridional $\phi = 0$ slices of the gyroscopic basis $\Phi_{m,l,k}^{\left(-\frac{1}{2},0\right)}$ for a parabolic cylinder, $h(s) = \frac{1}{5} \left( 1 + 4 s^{2} \right)$.
    The left half $(x < 0)$ of each subplot is the physical $z$ domain, while the right half $(x > 0)$ is the equivalent stretched $\eta$ domain.
    The index $(m,l,k)$ appears as text in each subplot.  Rows correspond to constant $(m,l)$ choices while columns correspond to constant $k$.
    }
\label{fig:cylinder_modes}
\end{figure}

In Figure~\ref{fig:cylinder_modes} we display some select gyroscopic modes for a cylinder with parabolic height.  The index $l$ denotes the vertical complexity while the index $k$ denotes the radial complexity.

\subsection{Stretched Annulus Bases}
We proceed analogously to the stretched cylinder case, using the appropriate volume element.
We weight the volume element $\dd V$ with the boundary polynomial 
\begin{equation}
m^{(\alpha)}(t, \eta) \triangleq \left( 1 - \eta^{2} \right)^{\alpha} \left( 1 - t^{2} \right)^{\alpha} \htilde(t)^{2 \alpha},
\label{eqn:malpha_annulus}
\end{equation}
where the only difference compared to the cylinder case in (\ref{eqn:malpha_cylinder}) is an additional $\left( 1 + t \right)^{\alpha}$ factor corresponding to the inner annulus boundary at $s = S_{i}$.
Define the weighted measure $\dd \mu(\alpha)$ for $\alpha > -1$ such that
\begin{equation}
\begin{aligned}
\dd \mu (\alpha) &\triangleq \frac{1}{2 \pi} \frac{4}{S_{o}^{2} - S_{i}^{2}} \left( 1-\eta^{2} \right)^{\alpha} \left( 1 - t^{2} \right)^{\alpha} \htilde(t)^{2 \alpha} \dd V \\
    &= \frac{1}{2 \pi} \left( 1-\eta^{2} \right)^{\alpha} \left( 1 - t^{2} \right)^{\alpha} \htilde(t)^{2 \alpha + 1 } \dd \phi \dd \eta \dd t
\end{aligned}
\end{equation}
and induced inner product $\langle \cdot , \cdot \rangle_{\dd \mu(\alpha)}$.

Though there are no radial coordinate singularities we still insert an $s^{m}$ prefactor in front of our basis
functions.  This guarantees the basis functions match the cylinder case as $S_{i} \to 0$ and importantly prevents
node clustering around $s = S_{i}$ for small $S_{i}$.  We note
\begin{equation}
2 s^{2} = S_{i}^{2} (1-t) + S_{o}^{2} (1+t)
\end{equation}
and define
\begin{equation}
\stilde(t) = \left( S_{i}^{2} (1-t) + S_{o}^{2} (1+t) \right)^{\frac{1}{2}}
\end{equation}
Now we define the gyroscopic basis as follows:
\begin{equation}
\begin{aligned}
\Upsilon_{m,l,k}^{(\alpha,\sigma)} &= e^{i m \phi} \stilde(t)^{\left|m\right|+\sigma} \htilde(t)^{l} P_{l}^{(\alpha,\alpha)}(\eta)
     Q_{k}^{\left( \alpha, \alpha, 2l + 2\alpha + 1, \left|m\right|+\sigma; \htilde, \stilde \right)}(t),
\label{eqn:annulus_basis}
\end{aligned}
\end{equation}
so that 
\begin{equation}
\left \langle \Upsilon_{m,l,k}^{(\alpha,\sigma)}, \Upsilon_{m',l',k'}^{(\alpha,\sigma)} \right \rangle_{\dd \mu(\alpha)} = \delta_{m,m'} \, \delta_{l,l'} \, \delta_{k,k'}.
\end{equation}

\begin{figure}[H]
\centering
\includegraphics[width=\linewidth]{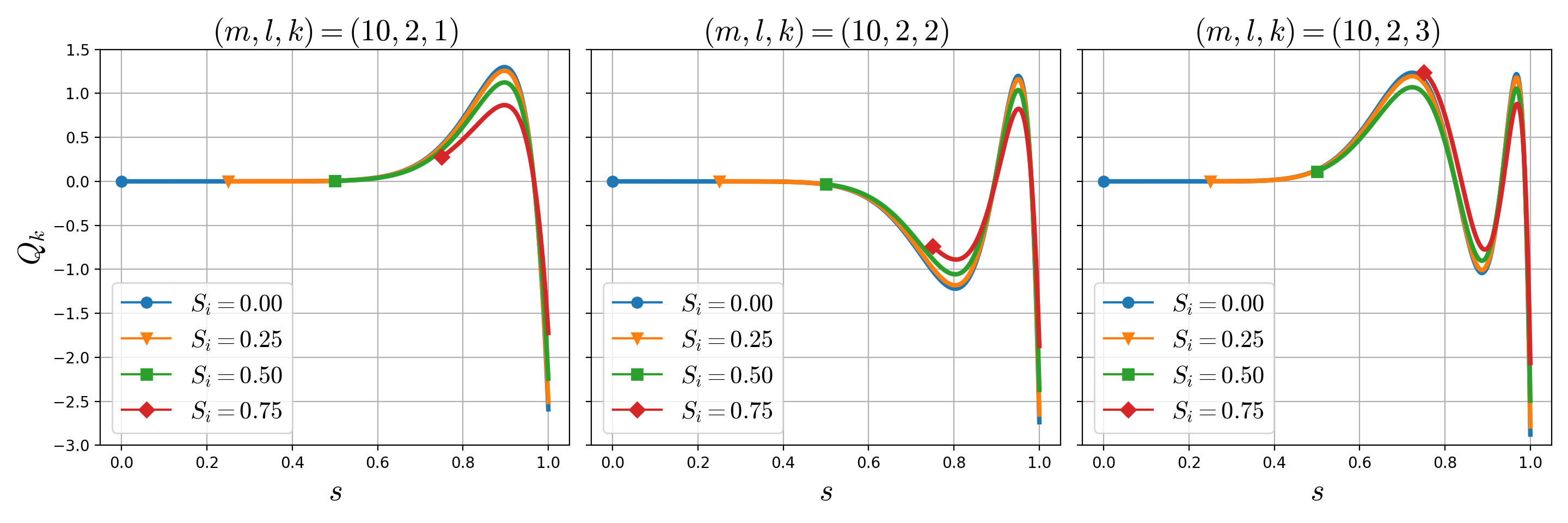}
\caption{Radial part of the cylinder and annulus basis functions for various $S_{i}$, $(\alpha, \sigma) = (0,0)$.  Observe that for small $S_{i}$ the bases
    approach those of the full cylinder $S_{i} = 0$.}
\label{fig:annulus_radial_part}
\end{figure}

For $\alpha = 0$ the annulus basis functions limit to the cylinder ones as $S_{i} \to 0$:
\begin{equation}
\lim_{S_{i} \to 0} \Upsilon_{m,l,k}^{(0,\sigma)} = \Phi_{m,l,k}^{(0,\sigma)}.
\end{equation}
We plot the radial generalized Jacobi polynomial in Figure~\ref{fig:annulus_radial_part}, where we see visual proof of the limiting case.
In particular, Figure~\ref{fig:annulus_radial_part} shows curves for azimuthal mode $m = 10$.  All curves decay like $s^{m} = s^{10}$ as $s \to 0$ \emph{independent of whether $s = 0$ is in the domain}.
This ensures the radial nodes are clustered well away from the $s = S_{i}$ boundary when $S_{i}$ is small, making the resolution more uniform throughout the domain.
Mahajan \cite{Mahajan} explored Zernike-type polynomials for the annulus with similar $S_{i} \to 0$ limiting properties.

\subsection{Half Domains}
\begin{figure}[H]
\centering
\includegraphics[width=\linewidth]{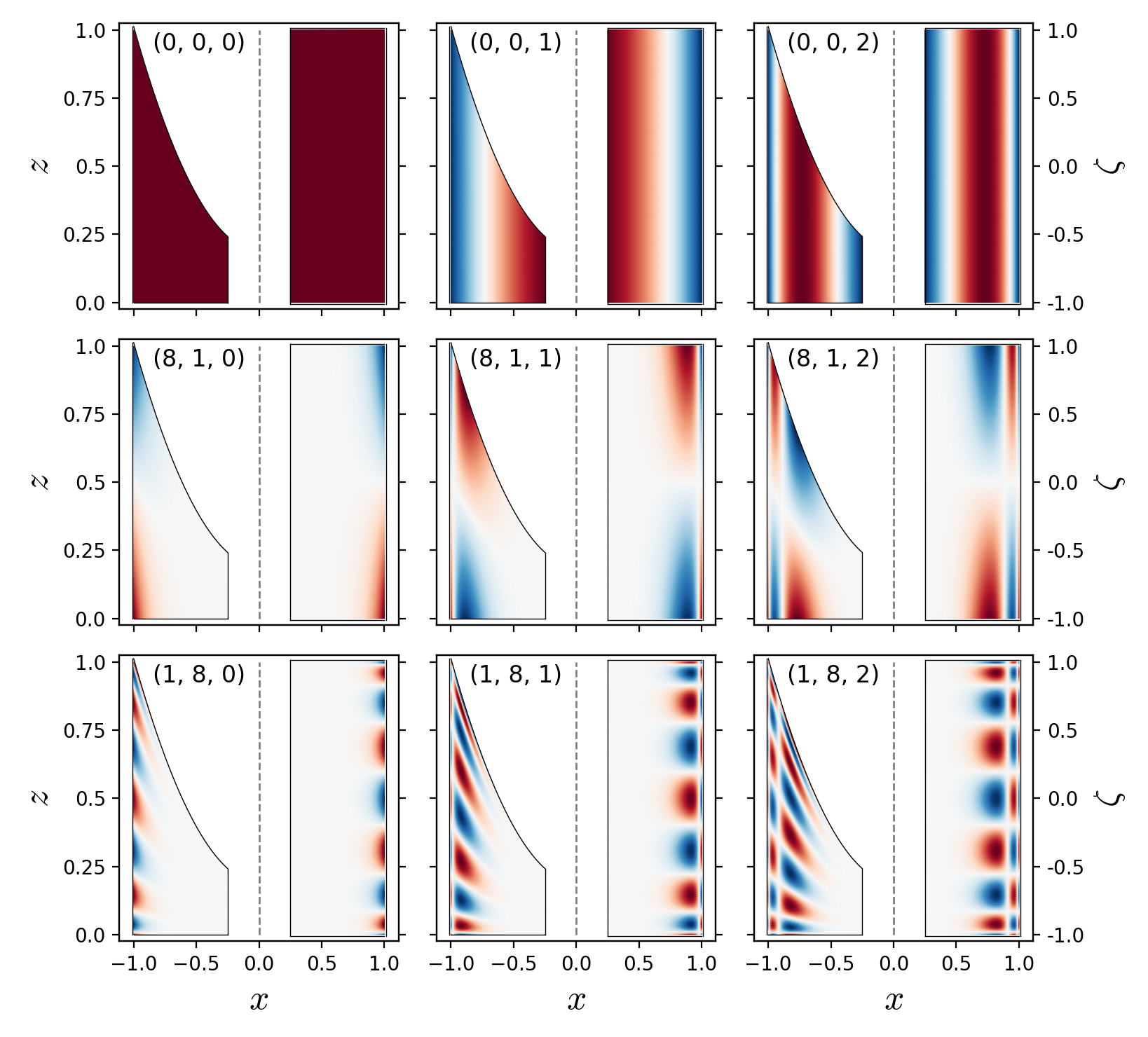}
\caption{Meridional $\phi = 0$ slices of the gyroscopic basis $\Upsilon_{m,l,k}^{\left(-\frac{1}{2},0\right)}$ for a parabolic annulus with radii $[0.25, 1.0]$.
    The left half $(x < 0)$ of each subplot is the physical $z$ domain, while the right half $(x > 0)$ is the stretched $\zeta$ domain.
    The index $(m,l,k)$ appears as text in each subplot.  Rows correspond to constant $(m,l)$ choices while columns correspond to constant $k$.
    }
\label{fig:annulus_modes}
\end{figure}

Using the $\zeta$ coordinate as in Section~\ref{sec:upper_half_geometries} we can implement an orthonormal basis hierarchy for the upper-half domains.
Recall that to use the $\zeta$ coordinate the domain must have no equatorial singularities and no square root around $\htilde$.
For both the full cylinder and annulus we simply replace the Jacobi polynomials in $\eta$ with polynomials in $\zeta$ while maintaining the rest of the
basis structure.  
In Figure~\ref{fig:annulus_modes} we plot a few upper-half domain annulus modes.

\subsection{Field Expansions}
We have now defined the hierarchy of gyroscopic bases to discretize scalar and vector fields in stretched cylindrical and annular domains.
As a reminder, we will expand our fields with basis functions $\Psi$, where $\Psi$ denotes either the cylinder basis $\Phi$ or the annulus basis $\Upsilon$.
Scalar fields decompose into their gyroscopic polynomial components:
\begin{equation}
f(t, \phi, \eta) = \sum_{m = -\infty}^{\infty}  \sum_{l = 0}^{\infty} \sum_{k = 0}^{\infty} \Psi_{m,l,k}^{(\alpha,0)}(t,\phi,\eta) \widehat{F}_{m,l,k}^{(\alpha)}
\label{eqn:scalar-expansion}
\end{equation}
where we freely specify $\alpha > -1$.
Vector fields decompose into a sum of their spin components:
\begin{equation}
\vec{u}(t,\phi,\eta) = \sum_{\sigma} \vec{\hat{e}}_{\sigma} \sum_{m = -\infty}^{\infty} \sum_{l = 0}^{\infty} \sum_{k = 0}^{\infty} \Psi_{m,l,k}^{(\alpha, \sigma)}(t, \phi, \eta) \widehat{U}_{m,l,k}^{(\alpha,\sigma)}
\end{equation}
The coefficient vectors $\widehat{F}_{m,l,k}^{(\alpha)}$ and $\widehat{U}_{m,l,k}^{(\alpha,\sigma)}$ are the expansion coefficients in the gyroscopic basis.
Using the generalized Jacobi polynomial algebra of Section~\ref{sec:generalized_jacobi} we obtain sparse matrix representations of calculus operators
acting on these coefficient vectors.

We compute the expansion coefficients using the inner product induced by $\dd \mu(\alpha)$:
\begin{equation}
\widehat{F}_{m,l,k}^{(\alpha)} = \left\langle f, \Psi_{m,l,k}^{(\alpha,0)} \right\rangle_{\dd \mu(\alpha)}
\end{equation}
and
\begin{equation}
\widehat{U}_{m,l,k}^{(\alpha,\sigma)} = \left\langle \vec{e}_{\sigma}^{*} \cdot \vec{u}, \Psi_{m,l,k}^{(\alpha,\sigma)} \right\rangle_{\dd \mu(\alpha)}.
\end{equation}
We approximate the integrals in the inner products using the Jacobi-$(\alpha,\alpha)$ quadrature rule in the vertical coordinate followed by the appropriate generalized Jacobi
quadrature rule in the radial coordinate.  This process utilizes bases' nested orthogonality to sift out the radial dependence attached to each azimuthal and and vertical mode.

\section{Discretization}\label{sec:discretization}
We now describe how to use the gyroscopic basis functions to discretize PDEs in the stretched cylindrical and annular domains.
This strategy is completely parallel to that developed for the Spherinder geometry \cite{Ellison_Julien_Vasil_2022}.
In fact, the gyroscopic bases for cylindrical domains is a direct generalization of the spherinder basis.
We see the gyroscopic basis reduce to the spherinder basis when we choose height function $h(t) = \sqrt{1-t}$.

\subsection{Operators}\label{sec:operators}
Generalized Jacobi polynomial algebra enables sparse matrix implementations of vector calculus.
Operators map between  the ${ (l,k) }$ indices of the basis functions as well as ${ (\sigma, \alpha) }$.
When an operator maps between scalars and vectors, like the scalar gradient or vector divergence, we see a natural
decomposition into the three components: a $\sigma$-raising operator, a $\sigma$-lowering operator and a $\sigma$-preserving operator.

\subsubsection{Regularity}
Paralleling \cite{Ellison_Julien_Vasil_2022}, we describe the regularity structure of fields using the \emph{regularity space} of degree $m$ defined by
\begin{equation}
\textrm{Reg}(m) = \left\{ f:[0,1] \to \mathbb{C} \hspace{2ex} \text{ s.t. } \hspace{2ex} 
    f(s) \sim s^{\left|m\right|} F\left(s^2\right) \hspace{1ex} \text{ as } \hspace{1ex} s \to 0
    \right\},
\end{equation}
where $F(s^2)$ is any even function of $s$ that is analytic in neighborhood of $s = 0$.
Then the degree $m$ azimuthal mode of a scalar field lives in $\textrm{Reg}(m)$ while vector fields decompose into 
the direct sum of three regularity spaces, $\textrm{Reg}(m - 1) \oplus \textrm{Reg}(m) \oplus \textrm{Reg}(m + 1)$
\cite{Vasil_Burns_Lecoanet_Olver_Brown_Oishi_2016}, corresponding to the spin components of the vector field.

We now define a hierarchy of Hilbert spaces indexed by real parameter $\alpha > -1$:
\begin{equation}
\mathcal{H}^{\alpha}(m) = \left\{ f \in \textrm{Reg}(m) \hspace{2ex}\text{ s.t. }\hspace{2ex} \norm{f}_{\dd \mu(\alpha)} < \infty \right\}.
\end{equation}
Differential operators map between these Hilbert spaces.  Generalized Jacobi polynomial algebra provides a sparse
representation of these operators acting on the gyroscopic bases.

These definitions exactly match those of \cite{Ellison_Julien_Vasil_2022}.
The regularity structure at $s = 0$ is identical in both cases and hence our definition of $\mathcal{H}^{\alpha}(m)$ matches up to the
definition of the norm $\norm{ \cdot }_{\dd \mu(\alpha)}$.

\subsubsection{Differential Operators}
In \cite{Ellison_Julien_Vasil_2022} we defined the gradient, divergence and curl operators in a coordinate-invariant way.
When we projected vector fields onto the the spinor basis we saw these differential operators split into three
fundamental parts: a spin-raising operator, a spin-lowering operator and a spin-preserving operator.  Because 
the generalized gyroscopic bases have the same regularity structure as in the spherinder we observe the same
underlying structure.

We now define the \emph{fundamental differential operators} $\mathcal{D}^{\delta}$ for $\delta \in \left\{ \pm 1, 0 \right\}$ such that
\begin{equation}
\mathcal{D}^{\delta}_{(\alpha,\sigma)} : \mathcal{H}^{\alpha}(m+\sigma) \mapsto \mathcal{H}^{\alpha+1}(m+\sigma+\delta),
\end{equation}
where
\begin{equation}
\begin{aligned}
\mathcal{D}^{\delta}_{(\alpha,\sigma)} &\triangleq \nabla_{\delta}
= \begin{cases} \frac{1}{\sqrt{2}} \left( \frac{\partial}{\partial S} \mp \frac{m}{S} \right) &\delta = \pm 1 \\[1ex]
    \frac{\partial}{\partial Z} &\delta = 0\end{cases}
\end{aligned}
\end{equation}
To compute the operator coefficients we need the constants $\gamma_{l}^{(\alpha)}$, $\delta_{l}^{(\alpha)}$ and $\lambda_{l}^{(\alpha)}$ such that
\begin{equation}
\begin{aligned}
P_{l}^{(\alpha,\alpha)}(\eta) &= \gamma_{l}^{(\alpha)} P_{l}^{(\alpha+1,\alpha+1)}(\eta) - \delta_{l}^{(\alpha)} P_{l-2}^{(\alpha+1,\alpha+1)}(\eta) \\
\frac{d}{d\eta} P_{l}^{(\alpha,\alpha)}(\eta) &= \lambda_{l}^{(\alpha)} P_{l-1}^{(\alpha+1,\alpha+1)}(\eta).
\label{eqn:jacobi_conversion_coefficients}
\end{aligned}
\end{equation}
The fundamental operators in the full cylinder then have the following structure:
\begin{equation}
\begin{aligned}
\frac{S_{o}}{2} \mathcal{D}^{+}_{(\alpha,\sigma)} \Phi_{m,l,k}^{(\alpha,\sigma)}
    &=  \phantom{-} \Phi_{m,l,\bigcdot}^{(\alpha+1,\sigma+1)}   \gamma_{l}^{(\alpha)} \bigg[ \mathcal{I}_{h} \, \mathcal{D}(+1,+1,+1) \bigg]_{\bigcdot,k} \\
    &\phantom{=} -  \Phi_{m,l-2,\bigcdot}^{(\alpha+1,\sigma+1)} \delta_{l}^{(\alpha)} \bigg[ \mathcal{I}_{h}^{\dagger} \, \mathcal{D}(+1,+1,-1) \bigg]_{\bigcdot,k} \\
\frac{S_{o}}{2} \mathcal{D}^{-}_{(\alpha,\sigma)} \Phi_{m,l,k}^{(\alpha,\sigma)}
    &=  \phantom{-} \Phi_{m,l,\bigcdot}^{(\alpha+1,\sigma-1)}   \gamma_{l}^{(\alpha)} \bigg[ \mathcal{I}_{h} \, \mathcal{D}(+1,-1,+1) \\
    &\phantom{=} -  \Phi_{m,l-2,\bigcdot}^{(\alpha+1,\sigma-1)} \delta_{l}^{(\alpha)} \bigg[ \mathcal{I}_{h}^{\dagger} \, \mathcal{D}(+1,-1,-1) \bigg]_{\bigcdot,k} \\
\mathcal{D}^{0}_{(\alpha,\sigma)} \Phi_{m,l,k}^{(\alpha,\sigma)}
    &=  \phantom{-} \Phi_{m,l-1,\bigcdot}^{(\alpha+1,\sigma)}   \lambda_{l}^{(\alpha)} \bigg[ \mathcal{I}_{a} \bigg]_{\bigcdot,k},
\end{aligned}
\end{equation}
where $\mathcal{I}_{a}$ is embedding with respect to the Jacobi $a$ parameter, $\mathcal{I}_{h}$ is embedding with respect to the height polynomial $\htilde$ defined in Section~\ref{sec:jacobi_embedding_operators},
and the three-argument $\mathcal{D}$ the generalized Jacobi differential operator defined in Section~\ref{sec:jacobi_differential_operators}.  Analogous to Einstein summation notation, we sum over repeated dots $\bigcdot$ in a product expression.

The annulus operators have a near-identical structure, changed only to account for the spin-weight polynomial proportional to $s$ moving to the final
position:
\begin{equation}
\begin{aligned}
\frac{S_{o}^{2} - S_{i}^{2}}{2} \mathcal{D}^{+}_{(\alpha,\sigma)} \Upsilon_{m,l,k}^{(\alpha,\sigma)}
    &=  \phantom{-} \Upsilon_{m,l,\bigcdot}^{(\alpha+1,\sigma+1)}   \, \gamma_{l}^{(\alpha)} \, \bigg[ \mathcal{I}_{h} \, \mathcal{D}(+1,+1,+1,+1) \bigg]_{\bigcdot,k} \\
    &\phantom{=} -  \Upsilon_{m,l-2,\bigcdot}^{(\alpha+1,\sigma+1)} \, \delta_{l}^{(\alpha)} \, \bigg[ \mathcal{I}_{h}^{\dagger} \, \mathcal{D}(+1,+1,-1,+1) \bigg]_{\bigcdot,k} \\
\frac{S_{o}^{2} - S_{i}^{2}}{2} \mathcal{D}^{-}_{(\alpha,\sigma)} \Upsilon_{m,l,k}^{(\alpha,\sigma)}
    &=  \phantom{-} \Upsilon_{m,l,\bigcdot}^{(\alpha+1,\sigma-1)}   \, \gamma_{l}^{(\alpha)} \, \bigg[ \mathcal{I}_{h} \, \mathcal{D}(+1,+1,+1,-1) \bigg]_{\bigcdot,k} \\
    &\phantom{=} -  \Upsilon_{m,l-2,\bigcdot}^{(\alpha+1,\sigma-1)} \, \delta_{l}^{(\alpha)} \, \bigg[ \mathcal{I}_{h}^{\dagger} \, \mathcal{D}(+1,+1,-1,-1) \bigg]_{\bigcdot,k} \\
\mathcal{D}^{0}_{(\alpha,\sigma)} \Upsilon_{m,l,k}^{(\alpha,\sigma)}
    &= \phantom{-} \Upsilon_{m,l-1,\bigcdot}^{(\alpha+1,\sigma)}   \, \lambda_{l}^{(\alpha)} \, \bigg[ \mathcal{I}_{a} \, \mathcal{I}_{b} \bigg]_{\bigcdot,k},
\end{aligned}
\end{equation}
where $\mathcal{I}_{b}$ is embedding with respect to the Jacobi $b$ parameter.

In the upper-half geometry the spin-preserving operators, corresponding to $z$-differentiation, are multiplied by a factor of two, 
a reflection of the $\zeta$ coordinate having half the physical range of the $\eta$ coordinate.
In addition, the spin-modifying operators map an $l \mapsto l-1$ component which for the cylinder takes the form
\begin{equation}
\begin{aligned}
\left\langle \mathcal{D}^{+}_{(\alpha,\sigma)} \Phi_{m,l,k}^{(\alpha,\sigma)}, \Phi_{m,l-1,j}^{(\alpha+1,\sigma+1)} \right\rangle_{\dd \mu(\alpha+1)}
    = - 2 \lambda_{l}^{(\alpha)} \, \bigg[ \mathcal{H}^{\prime} \, \mathcal{I}_{a} \, \mathcal{I}_{b} \bigg]_{j,k} \\
\left\langle \mathcal{D}^{-}_{(\alpha,\sigma)} \Phi_{m,l,k}^{(\alpha,\sigma)}, \Phi_{m,l-1,j}^{(\alpha+1,\sigma-1)} \right\rangle_{\dd \mu(\alpha+1)}
    = - 2 \lambda_{l}^{(\alpha)} \, \bigg[ \mathcal{H}^{\prime} \, \mathcal{I}_{a} \, \mathcal{I}_{b}^{\dagger} \bigg]_{j,k},
\end{aligned}
\end{equation}
where we define $\mathcal{H}^{\prime} = \frac{d}{dt} \htilde(t)$ as multiplication by the derivative of the height function.
For the annulus we have
\begin{equation}
\begin{aligned}
\left\langle \mathcal{D}^{+}_{(\alpha,\sigma)} \Upsilon_{m,l,k}^{(\alpha,\sigma)}, \Upsilon_{m,l-1,j}^{(\alpha+1,\sigma+1)} \right\rangle_{\dd \mu(\alpha+1)}
    = - 2 \lambda_{l}^{(\alpha)} \, \bigg[ \mathcal{H}^{\prime} \, \mathcal{I}_{a} \, \mathcal{I}_{b} \, \mathcal{I}_{s} \bigg]_{j,k} \\
\left\langle \mathcal{D}^{-}_{(\alpha,\sigma)} \Upsilon_{m,l,k}^{(\alpha,\sigma)}, \Upsilon_{m,l-1,j}^{(\alpha+1,\sigma-1)} \right\rangle_{\dd \mu(\alpha+1)}
    = - 2 \lambda_{l}^{(\alpha)} \, \bigg[ \mathcal{H}^{\prime} \, \mathcal{I}_{a} \, \mathcal{I}_{b} \, \mathcal{I}_{s}^{\dagger} \bigg]_{j,k},
\end{aligned}
\end{equation}
where $\mathcal{I}_{s}$ is embedding with respect to the $\stilde$ polynomial.

Just by matching spin weight we can easily construct traditional calculus operators acting on the basis functions from the 
fundamental differential operator $\mathcal{D}^{\delta}$.
We have the \emph{scalar gradient}
\begin{equation}
\begin{aligned}
\vec{\hat{e}}_{+}^{*} \cdot \del \Psi_{m,l,k}^{(\alpha,0)} &= \mathcal{D}^{+}_{(\alpha,0)} \, \Psi_{m,l,k}^{(\alpha,0)}\\
\vec{\hat{e}}_{-}^{*} \cdot \del \Psi_{m,l,k}^{(\alpha,0)} &= \mathcal{D}^{-}_{(\alpha,0)} \, \Psi_{m,l,k}^{(\alpha,0)}\\
\vec{\hat{e}}_{0}^{*} \cdot \del \Psi_{m,l,k}^{(\alpha,0)} &= \mathcal{D}^{0}_{(\alpha,0)} \, \Psi_{m,l,k}^{(\alpha,0)},
\end{aligned}
\end{equation}
and the \emph{divergence}
\begin{equation}
\begin{aligned}
\del \cdot \left( \vec{\hat{e}}_{+} \Psi_{m,l,k}^{(\alpha,+)} \right) &= \mathcal{D}^{-}_{(\alpha,+)} \, \Psi_{m,l,k}^{(\alpha,+)}\\
\del \cdot \left( \vec{\hat{e}}_{-} \Psi_{m,l,k}^{(\alpha,-)} \right) &= \mathcal{D}^{+}_{(\alpha,-)} \, \Psi_{m,l,k}^{(\alpha,-)}\\
\del \cdot \left( \vec{\hat{e}}_{0} \Psi_{m,l,k}^{(\alpha,0)} \right) &= \mathcal{D}^{0}_{(\alpha,0)} \, \Psi_{m,l,k}^{(\alpha,0)}.
\end{aligned}
\end{equation}
Note that we must keep careful track of $\alpha$ and $\sigma$ as we apply the $\mathcal{D}^{\delta}$ operators.  Each differential operator
raises $\alpha$ by one and modifies $\sigma$ according to $\delta$.  Hence the \emph{scalar Laplacian} $\mathcal{L}_{(\alpha)} : \mathcal{H}^{\alpha}(m) \mapsto \mathcal{H}^{\alpha+2}(m)$ takes the form
\begin{equation}
\begin{aligned}
\mathcal{L}_{(\alpha)} \Psi_{m,l,k}^{(\alpha,0)} &\triangleq \del^{2} \Psi_{m,l,k}^{(\alpha,0)} = \del \cdot \left( \del \Psi_{m,l,k}^{(\alpha,0)} \right) \\
    &= \left( \mathcal{D}^{-}_{(\alpha+1,+1)} \mathcal{D}^{+}_{(\alpha,0)} + \mathcal{D}^{+}_{(\alpha+1,-1)} \mathcal{D}^{-}_{(\alpha,0)} + \mathcal{D}^{0}_{(\alpha+1,0)} \mathcal{D}^{0}_{(\alpha,0)} \right) \Psi_{m,l,k}^{(\alpha,0)}.
\end{aligned}
\end{equation}

The fundamental differential operators carve the path to computing the vector curl.
We first compute the output spin components of the curl acting on each spin component of a general vector separately.
This reveals the spin coupling structure of the operator as well as scale factors for each coupling, here $\pm i$.
We then form the full operator by selecting the appropriate $\mathcal{D}^{\delta}$ for each spin coupling.
The end result is
\begin{equation}
\begin{aligned}
\del \times \left( \vec{\hat{e}}_{+} \Psi_{m,l,k}^{(\alpha,+)} \right) &= \phantom{-} i \, \vec{\hat{e}}_{+} \mathcal{D}^{0}_{(\alpha,+)} \Psi_{m,l,k}^{(\alpha,+)} - i \, \vec{\hat{e}}_{0} \mathcal{D}^{-}_{(\alpha,+)}  \Psi_{m,l,k}^{(\alpha,+)} \\
\del \times \left( \vec{\hat{e}}_{-} \Psi_{m,l,k}^{(\alpha,-)} \right) &= -i \, \vec{\hat{e}}_{0} \mathcal{D}^{+}_{(\alpha,-)} \Psi_{m,l,k}^{(\alpha,-)} + i \, \vec{\hat{e}}_{-} \mathcal{D}^{0}_{(\alpha,-)}  \Psi_{m,l,k}^{(\alpha,-)} \\
\del \times \left( \vec{\hat{e}}_{0} \Psi_{m,l,k}^{(\alpha,0)} \right) &= -i \, \vec{\hat{e}}_{+} \mathcal{D}^{+}_{(\alpha,0)} \Psi_{m,l,k}^{(\alpha,0)} + i \, \vec{\hat{e}}_{-} \mathcal{D}^{-}_{(\alpha,0)}  \Psi_{m,l,k}^{(\alpha,0)}.
\end{aligned}
\end{equation}
We define the \emph{vector Laplacian} $\mathcal{L}_{(\alpha,\sigma)} : \mathcal{H}^{\alpha}(m+\sigma) \mapsto \mathcal{H}^{\alpha+2}(m+\sigma)$ by
\begin{equation}
\mathcal{L}_{(\alpha,\sigma)} \Psi_{m,l,k}^{(\alpha,\sigma)} \triangleq \vec{\hat{e}}_{\sigma}^{*} \cdot \del^{2} \vec{\hat{e}}_{\sigma} \Psi_{m,l,k}^{(\alpha,\sigma)}.
\end{equation}
We specifically chose to represent vectors in the spinor basis so that this operation is diagonal in the spin components.
Using the gradient, divergence and curl operators we compute the $\mathcal{L}_{(\alpha,\sigma)}$ with
\begin{equation}
\del^{2} \vec{u} = \del \left( \del \cdot \vec{u} \right) - \del \times \del \times \vec{u}.
\end{equation}

\begin{figure}[H]
\centering
\includegraphics[width=0.8\linewidth]{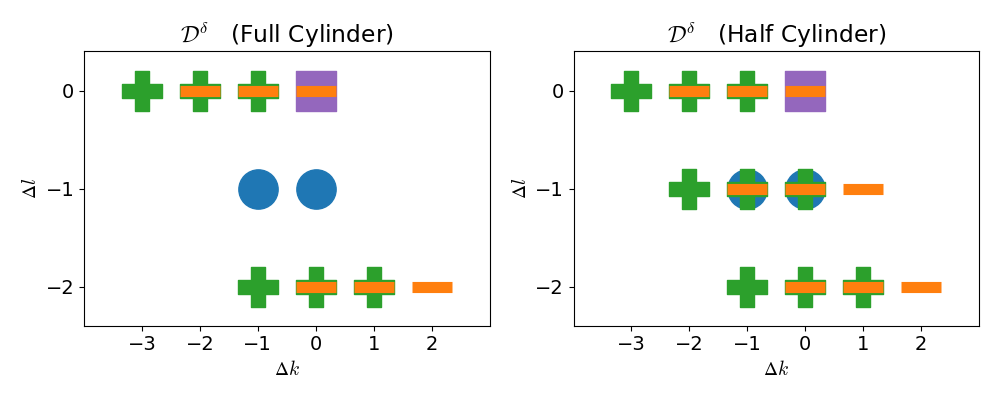}
\caption{Sparsity structure for the fundamental differential operators in cylindrical geometry with a parabolic height function.
The input mode is indicated by the square at location $(\Delta l, \Delta k) = (0,0)$.
The vertical axis $\Delta l$ is the \emph{change in vertical degree $l$} for the output of the operator, while the horizontal axis
$\Delta k$ is the \emph{change in horizontal degree $k$}.
Plus symbols indicate coupling for the $\delta = +1$ operator, minus symbols indicate coupling for $\delta = -1$ and disks indicate $\delta = 0$ coupling.
}
\label{fig:sparsity_diffops_cylinder}
\end{figure}

\begin{figure}[H]
\centering
\includegraphics[width=0.8\linewidth]{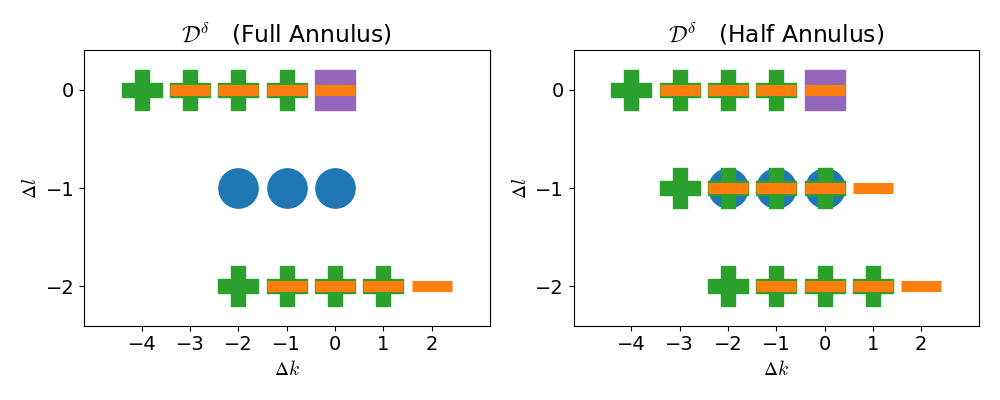}
\caption{Sparsity structure for the fundamental differential operators in annular geometry with a parabolic height function.
Bandwidth is one greater than that of the cylindrical domain of Figure~\ref{fig:sparsity_diffops_cylinder} despite matching height functions.
}
\label{fig:sparsity_diffops_annulus}
\end{figure}

We plot the sparsity structure of the three fundamental differential operators for the cylinder in Figure~\ref{fig:sparsity_diffops_cylinder} and
the annulus in Figure~\ref{fig:sparsity_diffops_annulus}.  The spin weight increment is given by the symbol in the figure - a plus indicates $\delta = +1$,
a minus indicates $\delta = -1$ and a circle indicates $\delta = 0$.  Observe the vertical $\zeta$ coordinate of the half domains generates an additional $l \mapsto l-1$ row for $\delta = \pm 1$
not present in the full geometry.  The annulus domain has bandwidth one greater than the cylinder domain due to its additional linear factor in the generalized Jacobi weight function for the radial
expansion.

\subsubsection{Coordinate Vector Operators}
Spin weight again guides operator definition for projection and coordinate vector multiplication just as for differential operators.
Multiplication by the radial $\vec{s} = s \, \vec{\hat{e}}_{S}$ vector takes the form
\begin{equation}
\begin{aligned}
\vec{s} \, \Phi_{m,l,k}^{(\alpha,0)}     &= \frac{S_{o}}{2} \vec{\hat{e}}_{+} \Phi_{m,l,\bigcdot}^{(\alpha,+)} \, \bigg[ \mathcal{I}_{b} \bigg]_{\bigcdot,k}
    + \frac{S_{o}}{2} \vec{\hat{e}}_{-} \Phi_{m,l,\bigcdot}^{(\alpha,-)} \,     \bigg[ \mathcal{I}_{b}^{\dagger} \bigg]_{\bigcdot,k} && \textrm{(Cylinder)} \\
\vec{s} \, \Upsilon_{m,l,k}^{(\alpha,0)} &= \frac{1}{2}     \vec{\hat{e}}_{+} \Upsilon_{m,l,\bigcdot}^{(\alpha,+)} \, \bigg[ \mathcal{I}_{s} \bigg]_{\bigcdot,k} 
    + \frac{1}{2}     \vec{\hat{e}}_{-} \Upsilon_{m,l,\bigcdot}^{(\alpha,-)} \, \bigg[ \mathcal{I}_{s}^{\dagger} \bigg]_{\bigcdot,k} && \textrm{(Annulus)},
\label{eqn:s_vector}
\end{aligned}
\end{equation}
showing the structure is identical when taking into account the spin weight factor is $1+t$ for cylinders and $\stilde$ for annuli.

Multiplication by the axial $\vec{z} = z \, \vec{\hat{e}}_{Z}$ vector produces a vector in the $\vec{\hat{e}}_{0}$ direction.
Recall that $\beta_{l}^{(\alpha)}$ is the off-diagonal term in the three-term recurrence for Jacobi-$(\alpha,\alpha)$ polynomials.
Then
\begin{equation}
\vec{z} \, \Psi_{m,l,k}^{(\alpha,0)} = \vec{\hat{e}}_{0} \left ( \beta_{l}^{(\alpha)} \Psi_{m,l+1,\bigcdot}^{(\alpha,0)} \, \bigg[ \mathcal{I}_{h}^{2} \bigg]_{\bigcdot,k} + \beta_{l-1}^{(\alpha)} \Psi_{m,l-1,\bigcdot}^{(\alpha,0)} \, \bigg[ \mathcal{I}_{h}^{\dagger \, 2} \bigg]_{\bigcdot,k} \right).
\label{eqn:z_vector_full}
\end{equation}
With the upper-half geometry we instead have
\begin{equation}
\vec{z} \, \Psi_{m,l,k}^{(\alpha,0)} = \frac{1}{2} \vec{\hat{e}}_{0} \left( \beta_{l}^{(\alpha)} \Psi_{m,l+1,\bigcdot}^{(\alpha,0)} \, \bigg[ \mathcal{I}_{h}^{2} \bigg]_{\bigcdot,k} 
    + \beta_{l-1}^{(\alpha)} \Psi_{m,l-1,\bigcdot}^{(\alpha,0)} \, \bigg[ \mathcal{I}_{h}^{\dagger \, 2} \bigg]_{\bigcdot,k}
    + \Psi_{m,l,\bigcdot}^{(\alpha,0)} \, \bigg[ \mathcal{I}_{h}^{\dagger} \, \mathcal{I}_{h} \bigg]_{\bigcdot,k} \right),
\end{equation}
which is identical to (\ref{eqn:z_vector_full}) except for the scale factor and $l \mapsto l$ term.

Dot products with the $\vec{s}$ and $\vec{z}$ vectors produce a spin-0 scalar field.  The operator structure is the same as (\ref{eqn:s_vector}) and (\ref{eqn:z_vector_full}) except
spin components move in the opposite directions.  That is, $\mathcal{I}_{b} \mapsto \mathcal{I}_{b}^{\dagger}$ in the cylinder and $\mathcal{I}_{s} \mapsto \mathcal{I}_{s}^{\dagger}$ in the annulus.

\begin{figure}[H]
\centering
\includegraphics[width=0.8\linewidth]{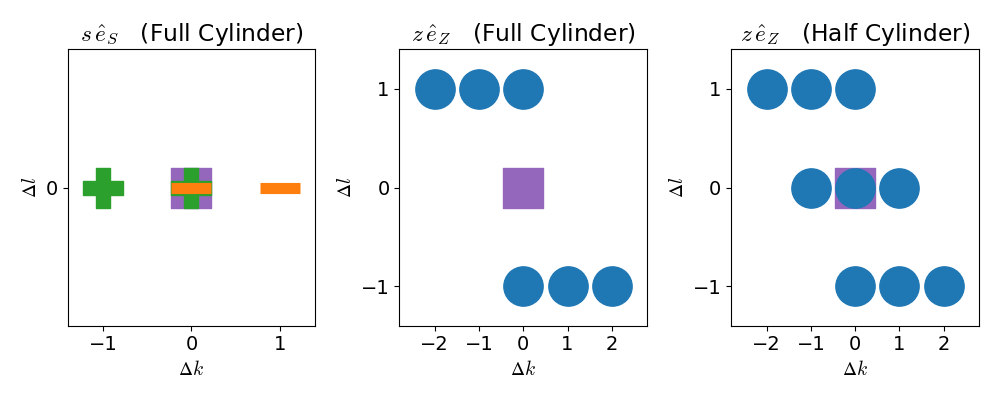}
\caption{Sparsity structure for the coordinate vector operators in cylindrical geometry with a parabolic height function}
\label{fig:sparsity_vecops}
\end{figure}

We plot the sparsity structure of the coordinate vector operators for the cylinder domain in Figure~\ref{fig:sparsity_vecops}.  The $\vec{s}$ multiplication operator does not
change the $l$ degree, while that $\vec{z}$ operator maps $l \mapsto l \pm 1$ in the full domain.  Again the vertical $\zeta$ coordinate of the half domain forces an additional
$l \mapsto l$ coupling not present in the full geometry.

\subsubsection{Conversion}\label{sec:gyropoly_conversion}
The last remaining tool we need is the $\alpha$-conversion operator.  We define 
\begin{equation}
\mathcal{I}_{(\alpha,\sigma)} : \mathcal{H}^{\alpha}(m+\sigma) \mapsto \mathcal{H}^{\alpha+1}(m+\sigma)
\end{equation}
via
\begin{equation}
\begin{aligned}
\mathcal{I}_{(\alpha,\sigma)} \Phi_{m,l,k}^{(\alpha,\sigma)} &\triangleq \Phi_{m,l,k}^{(\alpha,\sigma)} \\
    &= \Phi_{m,l,\bigcdot}^{(\alpha+1,\sigma)}   \, \gamma_{l}^{(\alpha)} \, \bigg[ \mathcal{I}_{a} \, \mathcal{I}_{h}^{2}             \bigg]_{\bigcdot,k}
     - \Phi_{m,l-2,\bigcdot}^{(\alpha+1,\sigma)} \, \delta_{l}^{(\alpha)} \, \bigg[ \mathcal{I}_{a} \, \mathcal{I}_{h}^{\dagger \, 2}  \bigg]_{\bigcdot,k} &&\textrm{(Cylinder)}  \\
\mathcal{I}_{(\alpha,\sigma)} \Upsilon_{m,l,k}^{(\alpha,\sigma)} &\triangleq \Upsilon_{m,l,k}^{(\alpha,\sigma)} \\
    &= \Upsilon_{m,l,\bigcdot}^{(\alpha+1,\sigma)}   \, \gamma_{l}^{(\alpha)} \, \bigg[ \mathcal{I}_{a} \, \mathcal{I}_{b} \, \mathcal{I}_{h}^{2}            \bigg]_{\bigcdot,k}
     - \Upsilon_{m,l-2,\bigcdot}^{(\alpha+1,\sigma)} \, \delta_{l}^{(\alpha)} \, \bigg[ \mathcal{I}_{a} \, \mathcal{I}_{b} \, \mathcal{I}_{h}^{\dagger \, 2} \bigg]_{\bigcdot,k} &&\textrm{(Annulus)}.
\end{aligned}
\end{equation}
We can now convert expansion coefficients for functions in $\mathcal{H}^{\alpha}(m)$ to the expansion coefficients for that same function in the $\mathcal{H}^{\alpha+1}(m)$ basis.
This is necessary when forming equations that equate derivatives to function values; we always need the equations to live in a single $\mathcal{H}^{\alpha}(m)$ space.

\begin{figure}[H]
\centering
\includegraphics[width=0.8\linewidth]{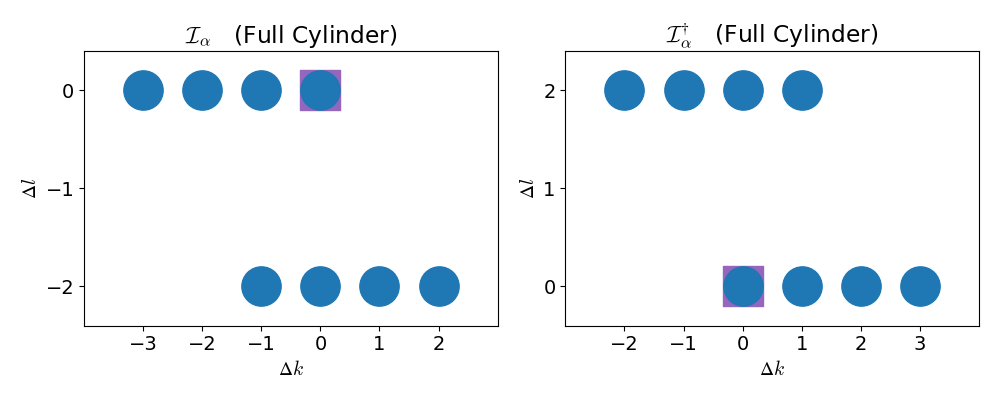}
\caption{Sparsity structure for the conversion operators in cylindrical geometry with a parabolic height function in $s$, which is linear in $t$.
Because the operators are adjoints of each other the sparsity structures are negated in both $\Delta l$ and $\Delta k$, i.e. $\Delta l \mapsto - \Delta l$ and $\Delta k \mapsto - \Delta k$.}
\label{fig:sparsity_conv}
\end{figure}

The embedding adjoint $\mathcal{I}_{(\alpha,\sigma)}^{\dagger} : \mathcal{H}^{\alpha}(m+\sigma) \mapsto \mathcal{H}^{\alpha-1}(m+\sigma)$
acts to multiply by the $\alpha$-weighted polynomial of the $\dd \mu(\alpha)$ measure:
\begin{equation}
\begin{aligned}
\mathcal{I}_{(\alpha,\sigma)}^{\dagger} \Phi_{m,l,k}^{(\alpha,\sigma)} &\triangleq \left(1-\eta^{2}\right) \left(1-t\right) \htilde(t)^{2} \Phi_{m,l,k}^{(\alpha,\sigma)} \\
    &= \Phi_{m,l,\bigcdot}^{(\alpha-1,\sigma)}   \, \gamma_{l}^{(\alpha)} \, \bigg[ \mathcal{I}_{a}^{\dagger} \, \mathcal{I}_{h}^{\dagger \, 2} \bigg]_{\bigcdot,k}
     - \Phi_{m,l+2,\bigcdot}^{(\alpha-1,\sigma)} \, \delta_{l}^{(\alpha)} \, \bigg[ \mathcal{I}_{a}^{\dagger} \, \mathcal{I}_{h}^{2}            \bigg]_{\bigcdot,k} &&\textrm{(Cylinder),}  \\
\mathcal{I}_{(\alpha,\sigma)}^{\dagger} \Upsilon_{m,l,k}^{(\alpha,\sigma)} &\triangleq \left(1-\eta^{2}\right) \left(1-t^{2}\right) \htilde(t)^{2} \Upsilon_{m,l,k}^{(\alpha,\sigma)} \\
    &= \Upsilon_{m,l,\bigcdot}^{(\alpha-1,\sigma)}   \, \gamma_{l}^{(\alpha)} \, \bigg[ \mathcal{I}_{a}^{\dagger} \, \mathcal{I}_{b}^{\dagger} \, \mathcal{I}_{h}^{\dagger \, 2} \bigg]_{\bigcdot,k}
     - \Upsilon_{m,l+2,\bigcdot}^{(\alpha-1,\sigma)} \, \delta_{l}^{(\alpha)} \, \bigg[ \mathcal{I}_{a}^{\dagger} \, \mathcal{I}_{b}^{\dagger} \, \mathcal{I}_{h}^{2}            \bigg]_{\bigcdot,k} &&\textrm{(Annulus)}.
\end{aligned}
\end{equation}
Here we see the maximum vertical degree increase $l \mapsto l+2$ as appropriate for multiplication by $1-\eta^{2}$.
We demonstrate the sparsity structure of the conversion operators for the cylinder domain in Figure~\ref{fig:sparsity_conv}.

\subsubsection{Triangular Truncation}
For the conversion operator $\mathcal{I}_{\alpha}$ to be an exact identity operator we need to triangularly truncate the radial expansion as a function of vertical degree.
Observe in the sparsity diagram for linear $\htilde$ the $l \mapsto l-2$ component has a $k \mapsto k+2$ piece.  This defines our radial truncation:
\begin{equation}
N_{\textrm{max}}(l) = N_{\textrm{max}}(0) - l \times \textrm{degree}\left(\htilde\right)
\end{equation}
with $N_{\textrm{max}}(0) \triangleq N_{\textrm{max}}$.  Then since $N_{\textrm{max}}(l-2) = N_{\textrm{max}}(l) + 2$ the operator is closed on our truncated expansion.
We employ $l$-dependent triangular truncation for all expansions using this basis, where truncation becomes steeper with increasing degree of $\htilde$.

\subsection{Boundary Evaluation}
\begin{figure}[H]
\centering
\includegraphics[width=\linewidth]{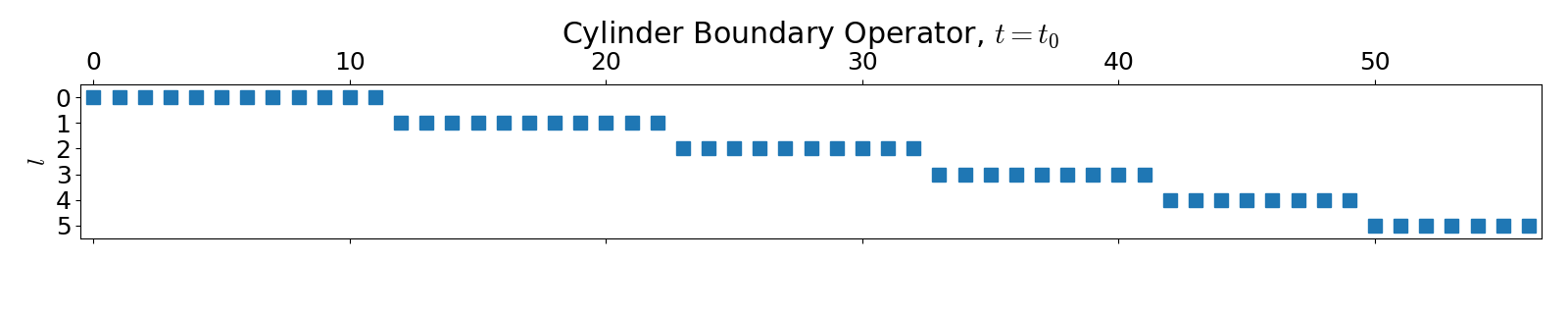}
\caption{Sparsity structure for boundary evaluation at $t = t_{0}$ with $(L_{\textrm{max}}, N_{\textrm{max}}) = (5, 11)$.  The row index corresponds to vertical degree $l$.}
\label{fig:boundary_s}
\end{figure}

Evaluating a field at a constant radius $t = t_{0}$ is straightforward.  We simply sum the radial series in $k$ evaluated at $t_{0}$ for each $l$.
For a field to vanish identically at $t = t_{0}$ we must have, for each $m$ and $l$,
\begin{equation}
    \sum_{k = 0}^{N_{\textrm{max}}(l)} Q_{k}^{\left(\cdots\right)}(t_0) \widehat{F}_{m,l,k}^{(\alpha,\sigma)} = 0
\end{equation}
We plot the sparsity structure of the $t = t_{0}$ boundary evaluation operator in Figure~\ref{fig:boundary_s}. 
We sort coefficients $\widehat{F}_{m,l,k}^{(\alpha,\sigma)}$ so that $k$ is the fastest index, then $l$, then $m$.
Notice the structure is block diagonal - the radial coefficients are weighted and summed for each $l$, corresponding to the rows of the operator. 

\begin{figure}[H]
\centering
\includegraphics[width=\linewidth]{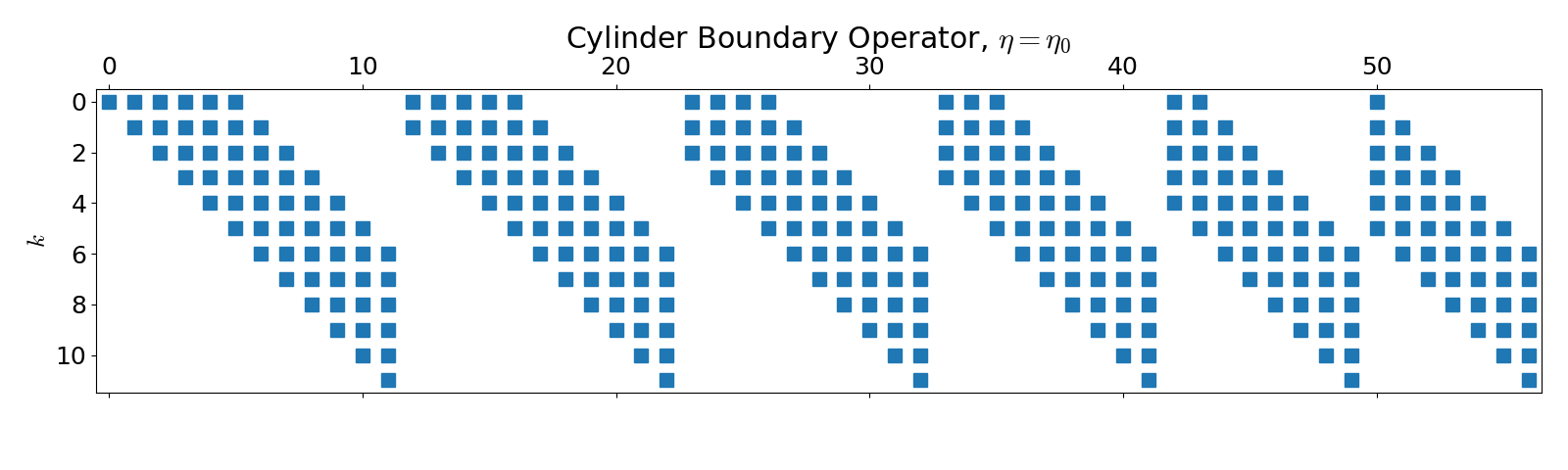}
\caption{Sparsity structure for boundary evaluation at $\eta = \eta_{0}$ (i.e., $z = h(s)$) with $(L_{\textrm{max}}, N_{\textrm{max}}) = (5, 11)$.  The row index corresponds to radial degree $k$
    while each block corresponds to increasing $l$.}
\label{fig:boundary_eta}
\end{figure}

Evaluating at a constant stretched height $\eta = \eta_{0}$ takes a bit more work.  We want to sum the vertical series in $l$ evaluated at $\eta_{0}$ for each $k$,
but the radial generalized Jacobi parameters depend on $l$.  To handle this we reverse the sums by converting the $l$-dependent $\htilde$ parameter to the common value $L_{\textrm{max}}+2\alpha+1$
using combinations of $\mathcal{I}_{h}$ and $\mathcal{I}_{h}^{\dagger}$.  Then to enforce a field vanishes at $\eta = \eta_{0}$ we have for each $m$ and $k$
\begin{equation}
\sum_{l = 0}^{L_{\textrm{max}}} P_{l}^{(\alpha,\alpha)}(\eta_{0}) \sum_{n=0}^{N_{\textrm{max}}(l)} \left[ \left( \mathcal{I}_{h} \right)^{L_{\textrm{max}}-l} \, \left( \mathcal{I}_{h}^{\dagger} \right)^{l} \right]_{k,n} \widehat{F}_{m,l,n}^{(\alpha,\sigma)} = 0.
\end{equation}
In contrast to the $t = t_{0}$ operator we see dense fill-in in the sparsity structure for evaluation at $\eta = \eta_{0}$, shown in Figure~\ref{fig:boundary_eta}.
In this case each row corresponds to a different $k$ index, with summation over $l$ that forces dense coupling of the basis functions.

\section{The Damped Inertial Waves Test Problem}\label{sec:eigenproblems}
We now demonstrate the power of the gyroscopic bases by exploring the damped inertial waves eigenproblem.
We choose a few different geometries that describe the Coreaboloid \cite{Lonner_Aggarwal_Aurnou_2022} at various rotation rates.
The upper equipotential surface of the Coreaboloid forms a paraboloid with steepness and depth parameterized
by the rotation rate and volume of the fluid.
We match the initial Coreaboloid experiment rotation rates and apply the gyroscopic discretization method to both cylindrical and annular domains.

\subsection{Damped Inertial Waves Discretization}
The damped inertial wave equations model the exponential decay in time of inertial modes for rotating fluids with viscosity.  
This system is a stepping stone toward the full convection problem which is thermally forced at the boundaries.
We solve the nondimensional Navier-Stokes equations (NSE)
\begin{equation}
\begin{aligned}
\partial_t \vec{u} + 2 \vec{\hat{e}}_z \vec{\times} \vec{u} &= -\del p + \textrm{E } \del^2 \vec{u} \\
\del \cdot \vec{u} &= 0
\end{aligned}
\end{equation}
in a rotating coordinate frame subject to no-slip boundary conditions ${ \vec{u} = 0 }$ at all bounding surfaces.  
Smallness of the Ekman number $\textrm{E} \ll 1$ indicates rotational forces dominate over viscous forces.
For convenience the Ekman number is defined as $\textrm{E} = \frac{\nu}{2 \Omega S_{o}^{2}}$ where $\nu$ denotes the fluid
viscosity and $\Omega$ the rotation rate.  Length, time and velocity are rescaled by $S_{o}$, $\Omega$, and $\Omega S_{o}$ respectively.
We substitute ${ \partial_t \mapsto \lambda }$ and solve for the generalized eigenvalues $\lambda$ of the discretized system.

To enforce the no-slip boundary condition we multiply the velocity field components by the boundary polynomial $m^{(1)}(t, \eta)$ defined in (\ref{eqn:malpha_cylinder}) for the cylinder and (\ref{eqn:malpha_annulus}) for the annulus.
This is a form of Galerkin recombination of the basis functions and is equivalent to the $\mathcal{I}_{(\alpha,\sigma)}^{\dagger}$ operator from Section~\ref{sec:gyropoly_conversion}.
Unsurprisingly, the boundary polynomial $m^{(1)}(t, \eta)$ vanishes on the boundary of the domain $t = \pm 1$, $\eta = \pm 1$, so fields recombined in this way explicitly satisfy the boundary condition.
We therefore define our velocity coefficient vectors $\widehat{U}^{(0,\sigma)}$ in terms of the auxiliary variable $\widehat{V}^{(1,\sigma)}$ via
\begin{equation}
\widehat{U}^{(0,\sigma)} \triangleq \mathcal{I}_{(1,\sigma)}^{\dagger} \widehat{V}^{(1,\sigma)}.
\end{equation}

The generalized eigensystem takes the form
\begin{equation}
\begin{pmatrix} L_{\textrm{bulk}} & L_{\textrm{proj}} \end{pmatrix} \begin{pmatrix} X \\ \widehat{\tau} \end{pmatrix} = \lambda \begin{pmatrix} M_{\textrm{bulk}} & M_{\textrm{proj}} \end{pmatrix} \begin{pmatrix} X \\ \widehat{\tau} \end{pmatrix}
\end{equation}
where
\begin{equation}
X =  \begin{pmatrix}[1.5] \mathcal{I}_{(1,+)}^{\dagger} \widehat{V}^{(1,+)} \\ \mathcal{I}_{(1,-)}^{\dagger} \widehat{V}^{(1,-)} \\ \mathcal{I}_{(1,0)}^{\dagger} \widehat{V}^{(1,0)} \\ \widehat{P}^{(1)} \end{pmatrix}, \hspace{8ex}
M_{\textrm{bulk}} = \begin{pmatrix}[1.2]
    \mathcal{I}_{(1,+)} \mathcal{I}_{(0,+)} & 0 & 0 & 0 \\ 
    0 & \mathcal{I}_{(1,-)} \mathcal{I}_{(0,-)} & 0 & 0 \\ 
    0 & 0 & \mathcal{I}_{(1,0)} \mathcal{I}_{(0,0)} & 0 \\ 
    0 & 0 & 0 & 0
\end{pmatrix},
\end{equation}
\begin{equation}
L_{\textrm{bulk}} = \begin{pmatrix}[1.5] 
    \textrm{E} \, \mathcal{L}_{(0,+)} - 2 i \, \mathcal{I}_{(1,+)} \mathcal{I}_{(0,+)} & 0 & 0 & -\mathcal{D}_{(1,0)}^{+}  \\ 
     0 & \textrm{E} \, \mathcal{L}_{(0,-)} + 2 i \, \mathcal{I}_{(1,-)} \mathcal{I}_{(0,-)} & 0 & -\mathcal{D}_{(1,0)}^{-} \\ 
     0 & 0 & \textrm{E} \, \mathcal{L}_{(0,0)} & -\mathcal{D}_{(1,0)}^{0}                                                  \\ 
     \mathcal{D}_{(0,+)}^{-} & \mathcal{D}_{(0,-)}^{+} & \mathcal{D}_{(0,0)}^{0} & 0                                      
\end{pmatrix},
\end{equation}
and
\begin{equation}
L_{\textrm{proj}} = 
    \textrm{diag} \left( 
        \mathcal{P}_{(2,+)},
        \mathcal{P}_{(2,-)},
        \mathcal{P}_{(2,0)},
        \mathcal{P}_{(1,0)}
        \right),
\hspace{8ex}
M_{\textrm{proj}} = \vec{0}.
\end{equation}
We solve this generalized eigenproblem for unknowns $\left( \widehat{V}^{(1,\sigma)}, \widehat{P}^{(1)}, \widehat{\tau} \right)$, where $\widehat{\tau}$ are the Lanczos-tau coefficients;
see Section~\ref{sec:damped_iwaves_projections} for an introduction to these additional variables and corresponding projection operators $\mathcal{P}_{(\alpha,\sigma)}$.

The momentum equations live in $\mathcal{H}^{2}(|m|+\sigma)$ while the divergence equations live in $\mathcal{H}^{1}(|m|)$.
The cascaded conversion operators in the $M$ matrix and in the Coriolis terms make the momentum equations consistent with $\alpha = 2$.
We choose $\alpha = 1$ for both the velocity and pressure field.  The pressure gradient then lives in $\alpha = 2$ along with the rest of the momentum equation, so no conversions are necessary.

\subsubsection{Projection Operators}\label{sec:damped_iwaves_projections}
Each $\mathcal{P}_{(\alpha,\sigma)}$ matrix projects the tau polynomials onto an equation living in $\mathcal{H}^\alpha(|m|+\sigma)$.  We define these projections as the highest radial and vertical modes of
the $\alpha$-conversion operator:
\begin{equation}
\mathcal{P}_{(\alpha,\sigma)} \triangleq \mathcal{I}_{(\alpha-1,\sigma)}[:],
\end{equation}
where $[:]$ denotes slicing the highest (most oscillatory) modes in the $s$ and $\eta$ directions.  The actual number of modes sliced depends on the domain - since the full cylinder only has a single radial boundary we take
only the final radial mode.  The annulus has two radial boundaries; this requires slicing the final two radial coefficients for each vertical degree.  In both domains we have boundaries at the top and bottom so we always slice the final two vertical modes.

Note that other choices of projections are possible.  For example we could define
\begin{equation}
\mathcal{P}_{(2,\sigma)} \triangleq \left( \mathcal{I}_{(1,\sigma)} \mathcal{I}_{(0,\sigma)} \right)[:]
\end{equation}
or even
\begin{equation}
\mathcal{P}_{(\alpha,\sigma)} \triangleq \mathcal{I}[:],
\end{equation}
where the no-subscript $\mathcal{I} : \mathcal{H}^{\alpha}(|m|+\sigma) \mapsto \mathcal{H}^{\alpha}(|m|+\sigma)$ is the identity operator.
These choices lead to non-trivial numerical differences in the solution and the optimal choice typically requires trial and error.

\subsection{Fundamental Mode vs. Rotation Rate in the Coreaboloid}
The Coreaboloid \cite{Lonner_Aggarwal_Aurnou_2022} refers to a rotating fluid layer with a free surface confined to an annular container.
The rotation rate determines the steepness of the parabolic equipotential upper surface and therefore the geometry of the problem.
We describe this free surface using parameters from the Coreaboloid paper:
\begin{equation}
h_{\text{core}}(s) = h_{0} + \frac{ \Omega^{2} }{ 2 g } s^{2},
\hspace{4ex}
h_{0} = H_{NR} - \frac{ \Omega^{2} S_{o}^{2} }{ 4 g }.
\end{equation}
The experimental apparatus has values $S_i = 10.2 \, \text{cm}$, $S_{o} = 37.25 \, \text{cm}$ and $H_{NR} = 17.08 \, \text{cm}$.
and $g = 9.81 \text{m}/\text{s}^2$.
To improve numerical conditioning of the algorithm we nondimensionalize height and radius by $S_o$ so that we have $s \in [0.27,1]$,
$H_{NR}/S_{o} = 0.46$ and
\begin{equation}
h(s) = \frac{1}{S_o} h_{\text{core}}(S_{o} \, s) = \frac{1}{S_o} \left( h_{0} + \frac{ \Omega^{2} S_{o}^{2} }{ 2 g } s^{2} \right).
\end{equation}

This rescaling keeps the height function bounded near one.  Large dynamic ranges of $h$ cause loss of precision when computing the
three-term recurrence of the gyroscopic polynomials.  For vertical degree $l$ the basis functions have a prefactor of $h(t)^{l}$.
This compounds largeness in $h$, leading to inaccuracy in the computed operators and solution.
Proper nondimensionalization - namely rescaling lengths with $S_{o}$ - avoids this issue.

\begin{figure}[H]
\centering
\includegraphics[width=0.7\linewidth]{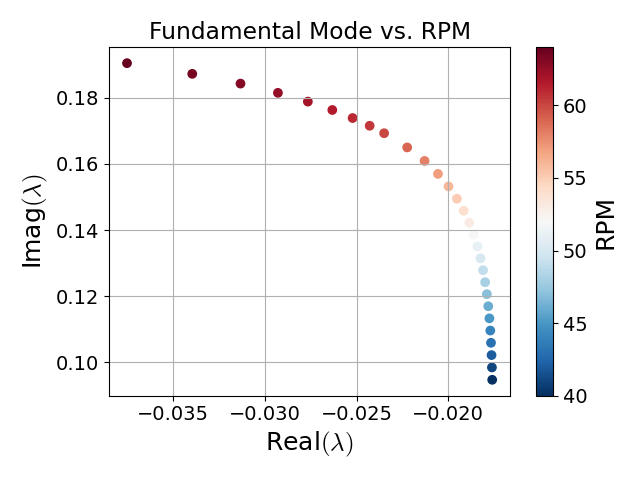}
\caption{Fundamental eigenvalue in the Coreaboloid for various rotation rates, measured in rotations per minute (RPM), with
$m = 14$ and $\textrm{E} = 10^{-5}$.
}
\label{fig:damped_iwaves-coreaboloid-rpm}
\end{figure}

\begin{figure}[H]
\centering
\includegraphics[width=\linewidth]{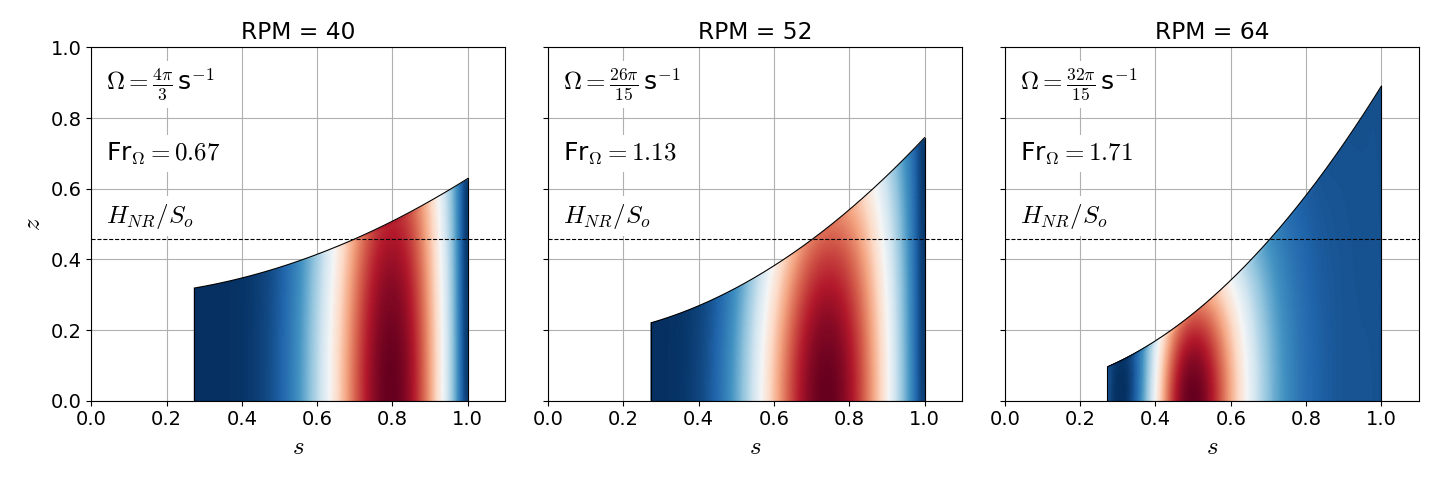}
\includegraphics[width=\linewidth]{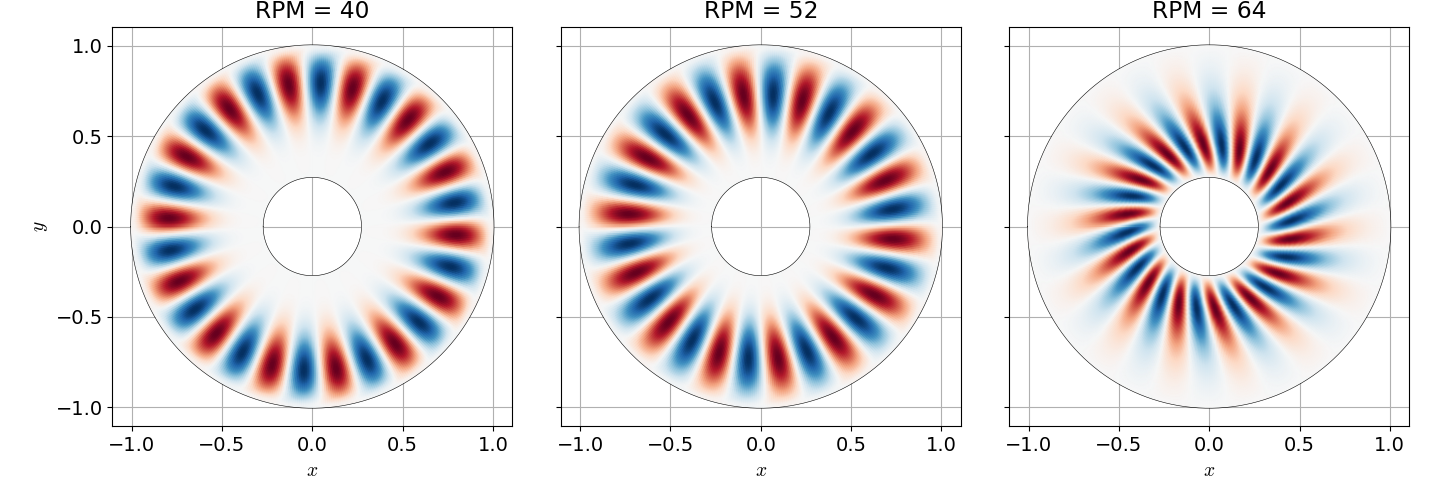}
\caption{Fundamental pressure mode in the Coreaboloid for $m = 14$ and $\textrm{E} = 10^{-5}$ at various rotation rates
viewed from the side, i.e. a meridional cross-section (top row), and from above (bottom row).
The geometry becomes singular above 66 RPM where the height at the inner radius vanishes.
We mark the nondimensionalized non-rotating height $H_{NR} / S_{o}$ that each surface intersects at $s = 1/\sqrt{2}$.
The rotational Froude number $Fr_{\Omega} = \frac{\Omega^{2} S_{o}}{g}$ measures the relative importance of the centrifugal force and gravity.
As the rotation rate increases the upper surface of the domain steepens.  We see the spiraling strengthen along with this increase in rate.
}
\label{fig:damped_iwaves-coreaboloid-rpm_p}
\end{figure}

In all experiments we set $m = 14$ and $\textrm{E} = 10^{-5}$.
Figure~\ref{fig:damped_iwaves-coreaboloid-rpm} shows the dependence of the critical mode on the rotation rate.
Figure~\ref{fig:damped_iwaves-coreaboloid-rpm_p} plots the fundamental pressure mode for select geometries.

\subsection{Fundamental Mode vs. Eccentricity in the Spheroid}
\begin{figure}[H]
\centering
\includegraphics[width=0.7\linewidth]{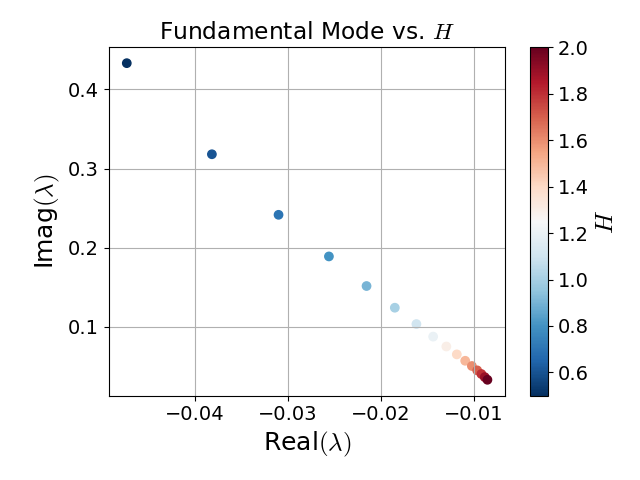}
\caption{Fundamental eigenvalue in the oblate spheroid for various heights $H$, where $m = 14$ and $\textrm{E} = 10^{-5}$.}
\label{fig:damped_iwaves-spheroid-height}
\end{figure}

\begin{figure}[H]
\centering
\includegraphics[width=\linewidth]{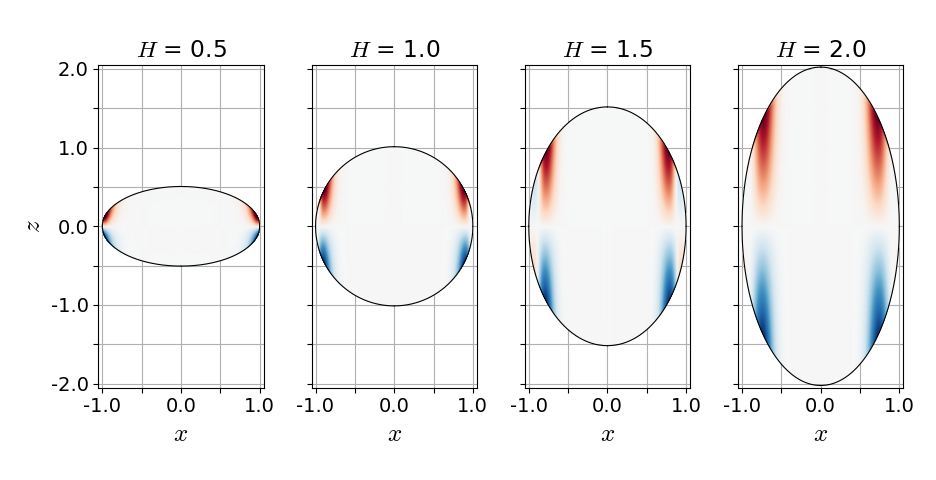}
\caption{Fundamental pressure mode in the oblate spheroid for select heights $H$.
The mode is localized at the equator for highly eccentric (low $H$) spheroids.  As $H$ increases the mode migrates into the bulk of the fluid.}
\label{fig:damped_iwaves-spheroid-height_p}
\end{figure}

We now compute the fundamental eigenmode of the damped inertial waves problem in an oblate spheroid.
We compute the solution for a range of eccentric spheres with height function
\begin{equation}
h(s) = H \, \sqrt{1 - \left( \frac{s}{S_{o}} \right)^{2}}.
\end{equation}
As above we set $m = 14$ and $\textrm{E} = 10^{-5}$ for all problems.
In Figure~\ref{fig:damped_iwaves-spheroid-height} we show the complex eigenfrequency $\lambda$ as a function of the height $H$.
Figure~\ref{fig:damped_iwaves-spheroid-height_p} demonstrates the pressure field of the mode at three select heights.

\subsection{Fundamental Mode vs. Inner Radius}
\begin{figure}[H]
\centering
\includegraphics[width=\linewidth]{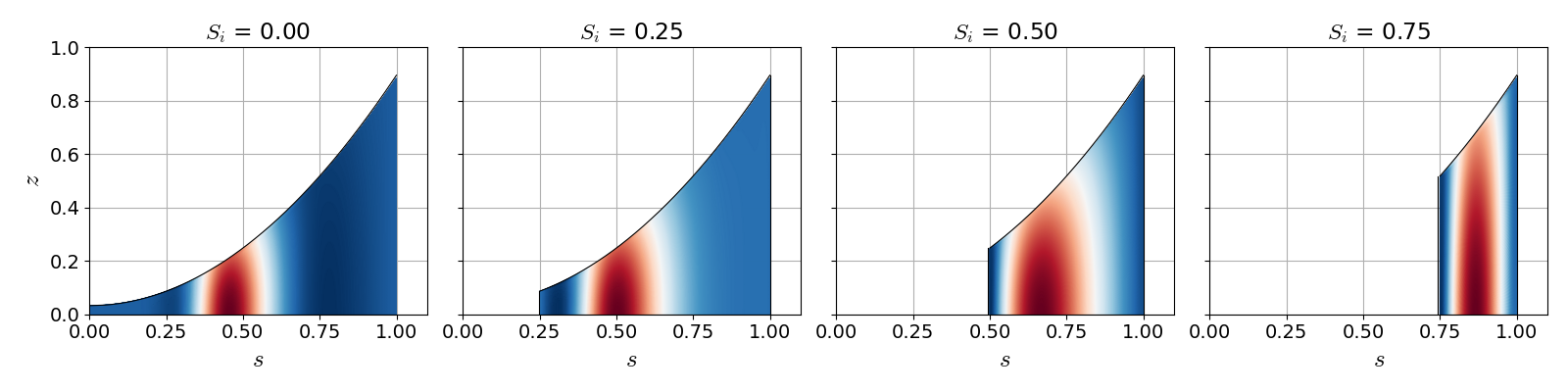}
\includegraphics[width=\linewidth]{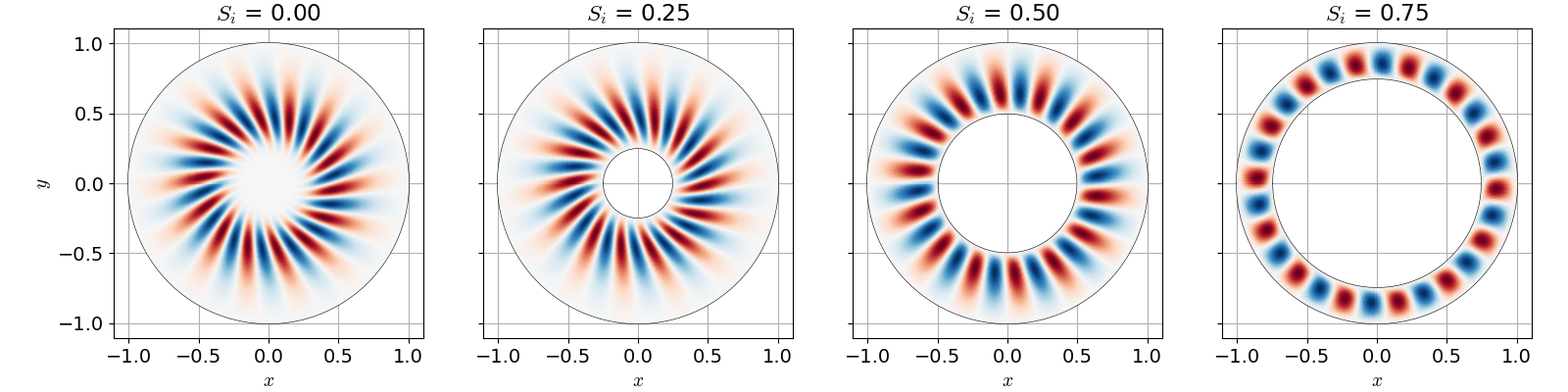}
\caption{Fundamental pressure mode in the Coreaboloid rotating at 64 RPM for select inner radii viewed from the side (top row) and from above (bottom row).
Spiraling is not evident in the modes with smaller inner radius $S_{i}$.
}
\label{fig:damped_iwaves-coreaboloid-radius_p}
\end{figure}

\begin{figure}[H]
\centering
\includegraphics[width=\linewidth]{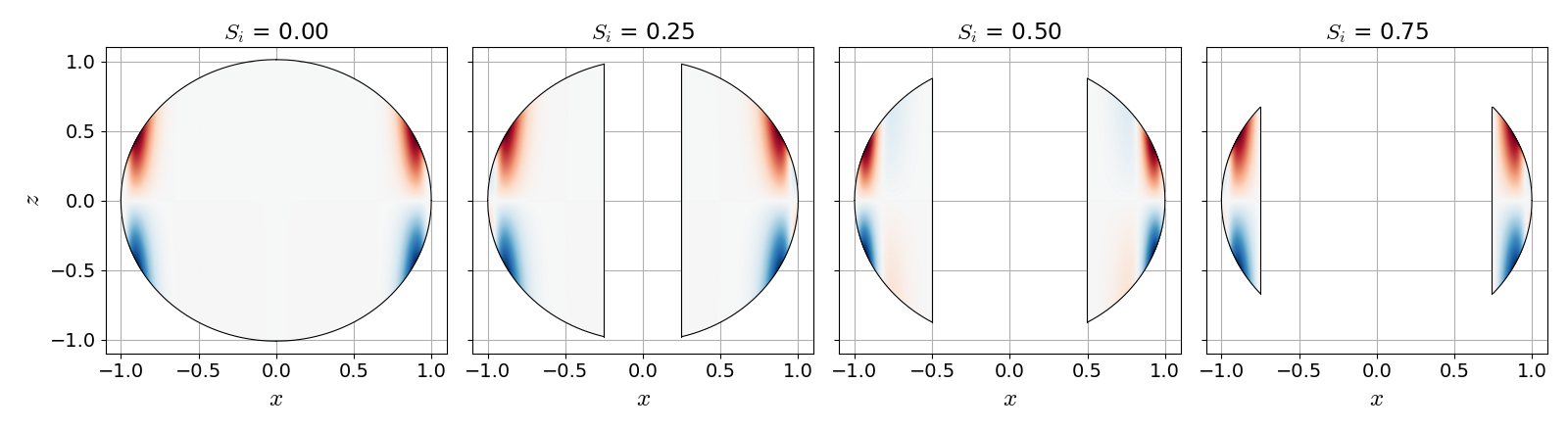}
\includegraphics[width=\linewidth]{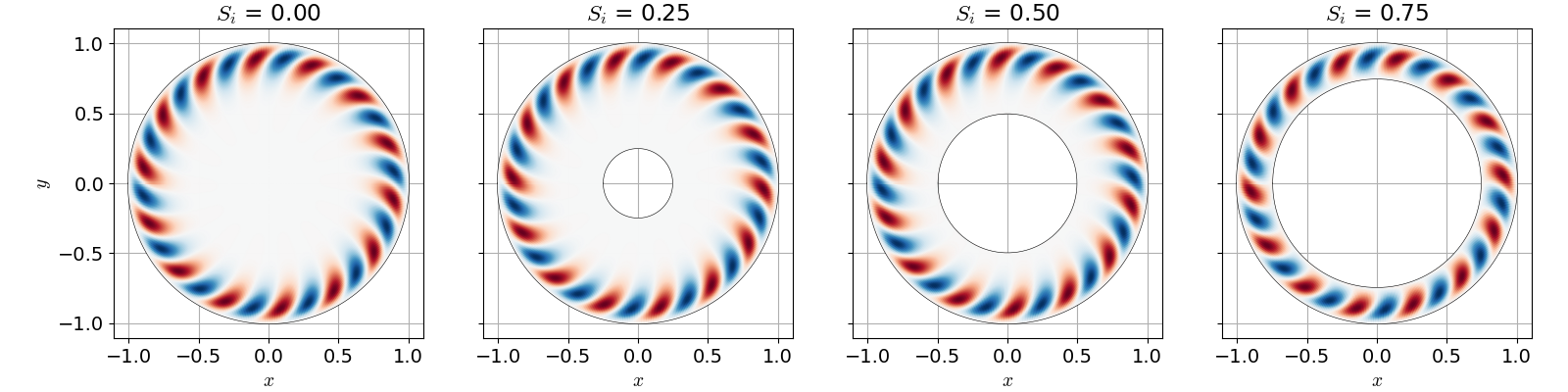}
\caption{Fundamental pressure mode in the sphere for select inner radii viewed from the side (top row) and from above (bottom row).  As in Figure~\ref{fig:damped_iwaves-coreaboloid-rpm_p} we see
evidence of spiraling.  In these cases the spiraling is \emph{retrograde} rather than prograde in the upper-half geometry.  This is because the upper-half pressure modes are
even in degree, while the the sphere modes are odd.  The different symmetries about $z = 0$ lead to the two different spiraling behaviors.
}
\label{fig:damped_iwaves-sphere-radius_p}
\end{figure}

We compute the fundamental mode for a fixed height profile as we vary the inner radius, holding constant both $m = 14$ and $\textrm{E} = 10^{-5}$.
We investigate both the Coreaboloid geometry at 64 RPM as well as the sphere with an inner cylinder along the $z$-axis excised.  In both cases
the fundamental mode is confined away from the central axis.  This means its frequency is relatively independent of the
inner radius.
Figure~\ref{fig:damped_iwaves-coreaboloid-radius_p} plots the fundamental pressure mode for select inner radii in the
Coreaboloid and Figure~\ref{fig:damped_iwaves-sphere-radius_p} plots it for the sphere.

\section{Conclusions}\label{sec:conclusion}
Motivated by the spatial anisotropy induced by rotating flows we introduce 
gyroscopic coordinate systems and corresponding orthogonal 3D bases for functional approximations in stretched cylinders and annuli.
The coordinate systems are tailored to gyroscopically aligned dynamics upon specification of the height function $z = h(s)$.
For polynomial-type $h(s)$ the method defines a natural hierarchy of bases that permit the sparse discretization of PDEs and hence
efficient, fast numerical computations.

The volume element of the gyroscopic coordinate system leads to a generalized Jacobi polynomial
weight function.  To proceed we defined generalized Jacobi polynomials with 
augmenting polynomial factors specified by $h(s)$.  Analogous to the Jacobi case we found sparse embedding and differential
operators that map polynomials from one set of generalized parameters to another, each parameter differing by an integer.
We developed efficient tools to calculate the polynomials and operator coefficients.  Importantly, the Christoffel-Darboux formulation for
computing the three-term recurrence relations provides a hierarchical means of calculating the polynomials.
In this way we compute all recurrences for a problem with fixed azimuthal order $m$ in $\mathcal{O}(N^{2})$ operations rather than
$\mathcal{O}(N^{3})$ using standard techniques.

With coordinate system and generalized Jacobi polynomials in hand we demonstrate how to construct the gyroscopic bases.
Because we used the volume element to define our basis the polynomials themselves are singularity free on their domain.
We chose particular coordinate-dependent prefactors so that the polynomials are truly Cartesian $(x,y,z)$ polynomials 
expressed in the gyroscopic coordinate system.
We found sparse differential operators mapping between bases in the hierarchy to implement tensor calculus for PDEs occurring in fluid dynamics.
We demonstrated their utility by solving the classical test eigenproblem of damped inertial waves in various domains.

The \textsc{Gyropoly} python package, available at \url{https://github.com/acellison/gyropoly}, contains tools for computing with generalized Jacobi polynomials
and scripts for reproducing all examples from this paper.  The python code provides a convenient interface for computing with generalized Jacobi polynomials.  On
top of this library we build the gyroscopic bases and associated calculus operators.  These tools are building blocks much in the same spirit as recent numerical
platforms such as the \textsc{Dedalus} software package \cite{Burns_Vasil_Oishi_Lecoanet_Brown_2020}.

\bigskip
\noindent \textbf{Funding: } This work was supported by the National Science Foundation Grant DMS-2009319.

\begin{appendices}

\section{Three-Term Recurrence for Quadratic Augmenting Factors}\label{app:christoffel_darboux}
Define
\begin{equation}
\begin{aligned}
w^{(a,b,c)}(x) &= (1-x)^{a} (1+x)^{b} \left(x^{2} + y^{2} \right)^c \\
w^{(a,b,c,d)}(x) &= (1-x)^{a} (1+x)^{b} (i y - x)^{c} (-i y - x)^{d}
\end{aligned}
\end{equation}
with
\begin{equation}
\begin{aligned}
\omega^{(a,b,c)}   &= \int_{-1}^{1}{ (1-x)^{a} (1+x)^{b} \left( x^{2} + y^{2} \right)^{c} \dd x} \\
\omega^{(a,b,c,d)} &= \int_{-1}^{1}{ (1-x)^{a} (1+x)^{b} (i y - x)^{c} (-i y - x)^{d} \dd x}
\end{aligned}
\end{equation}
Define the inner product
\begin{equation}
\left \langle p, q \right \rangle_{(a,b,c,d)} = \frac{1}{\omega^{(a,b,c,d)}} \int_{-1}^{1}{ w^{(a,b,c,d)}(x) \overline{ p(x) } q(x) \dd x},
\end{equation}
which generates orthonormal polynomials $P_{n}^{(a,b,c,d)}$ with three-term recurrence
\begin{equation}
x P_{n}^{(a,b,c,d)}(x) = \beta_{n}^{(a,b,c,d)} P_{n+1}^{(a,b,c,d)}(x) + \alpha_{n}^{(a,b,c,d)} P_{n}^{(a,b,c,d)}(x) + \beta_{n-1}^{(a,b,c,d)} P_{n-1}^{(a,b,c,d)}(x)
\end{equation}
We want to generate three-term recurrence coefficients $\alpha_{n}^{(a,b,c+1)} = \alpha_{n}^{(a,b,c+1,c+1)}, \beta_{n}^{(a,b,c+1)} = \beta_{n}^{(a,b,c+1,c+1)}$ from
known coefficients corresponding to $w^{(a,b,c)}$.  To do so we repeatedly apply the formulae for $\beta_{n}^{(a,b,c+1,d)}$, $\beta_{n}^{(a,b,c,d+1)}$, $\alpha_{n}^{(a,b,c+1,d)}$
and $\alpha_{n}^{(a,b,c,d+1)}$ in \cite{Snowball_Olver_2020a} to obtain
\begin{equation}
\beta_{n}^{(a,b,c+1)}
    = \left( \frac{ \left( \sum_{k=0}^{n} \left| P_{k}^{(a,b,c)}(i y) \right|^{2} \right) \left( \sum_{k=0}^{n+2} \left| P_{k}^{(a,b,c)}(i y) \right|^{2} \right) }{ \left( \sum_{k=0}^{n+1} \left| P_{k}^{(a,b,c)}(i y) \right|^{2} \right)^{2} } \right)^{\frac{1}{2}} \beta_{n+1}^{(a,b,c)}
\end{equation}
\begin{equation}
\begin{aligned}
\alpha_{n}^{(a,b,c+1)} 
    = \hspace{2ex}&\left( \frac{ \sum_{k=0}^{n+2} \left| P_{k}^{(a,b,c)}(i y) \right|^{2} }{ \sum_{k=0}^{n+1} \left| P_{k}^{(a,b,c)}(i y) \right|^{2} } - 1 \right) \textrm{Re} \left[ \frac{ P_{n+1}^{(a,b,c)}(i y) }{ P_{n+2}^{(a,b,c)}(i y) } \right] \beta_{n+1}^{(a,b,c)} \\
      - &\left( \frac{ \sum_{k=0}^{n+1} \left| P_{k}^{(a,b,c)}(i y) \right|^{2} }{ \sum_{k=0}^{n} \left| P_{k}^{(a,b,c)}(i y) \right|^{2} } - 1 \right) \textrm{Re} \left[ \frac{ P_{n}^{(a,b,c)}(i y)   }{ P_{n+1}^{(a,b,c)}(i y) } \right] \beta_{n}^{(a,b,c)}
      + \alpha_{n+1}^{(a,b,c)}.
\end{aligned}
\end{equation}
We bootstrap up to the desired $c$ parameter using the known recurrence and polynomials.  To generate the recurrence for polynomials of degree $N$ we require an initialization
with $N_{\textrm{init}} = N + 2(\textrm{floor}(c)+1)$ recurrence coefficients.  
When $c$ is an integer we generate the three-term recurrence for $P_{n}^{(a,b,c)}$ by starting with the $N_{\textrm{init}}$ terms in the Jacobi recurrence for $P_{n}^{(a,b)} \equiv P_{n}^{(a,b,0)}$.
When $c$ is not an integer we define $\gamma \triangleq c - \textrm{floor}(c)$ and compute $N_{\textrm{init}}$ terms in the recurrence for $P_{n}^{(a,b,\gamma)}$ using for example the Stieltjes procedure
or modified Chebyshev algorithm.  From this base stage we then use the Christoffel-Darboux formulations to recurse up to $c$.

This analysis works as well for arbitrary real quadratic factors that don't vanish in $(-1,1)$.  Simply factor the quadratic polynomial into its complex conjugate roots via
$x^{2} + c_{1} x + c_{2} = \left(x - z \right)\left(x - \overline{z}\right)$ and replace any occurrence of $i y$ above with $z$.  We therefore have a stable method of computing
the recurrence coefficients for arbitrary, non-vanishing polynomials with real coefficients.  Simply factor the polynomial into its roots and perform
the standard parameter increment on real zeros and the double-increment on conjugate pairs.

\section{Fully General Domain Formulations}\label{app:general_domains}

\subsection{The General Cylinder}
Define the weighted volume measure $\dd \mu(\alpha)$ for $\alpha > -1$ by
\begin{equation}
\begin{aligned}
\dd \mu (\alpha) &\triangleq \frac{1}{2 \pi} \frac{4}{S_{o}^{2}} \left( 1-\eta^{2} \right)^{\alpha} \left( 1 - t \right)^{\alpha} \htilde(t)^{2 \chi_{h} \alpha} \dd V \\
    &= \frac{1}{2 \pi} \left( 1-\eta^{2} \right)^{\alpha} \left( 1 - t \right)^{\frac{1}{2} \chi_{o} + \alpha} \htilde(t)^{\left( 2 \alpha + 1 \right) \chi_{h}} \dd \phi \dd \eta \dd t
\end{aligned}
\end{equation}
and inner product
\begin{equation}
\left \langle f, g \right \rangle_{\dd \mu(\alpha)} \triangleq \int \dd \mu(\alpha) \overline{f} g.
\end{equation}
We therefore define the basis function as follows:
\begin{equation}
\Phi_{m,l,k}^{(\alpha,\sigma)} = e^{i m \phi} (1+t)^{\frac{\left|m\right| +\sigma}{2}} (1-t)^{\frac{l}{2} \chi_{o}} \htilde(t)^{l \chi_{h}} P_{l}^{(\alpha,\alpha)}(\eta) \, Q_{k}^{\left( \left( l+\frac{1}{2} \right) \chi_{o} + \alpha,
\left|m\right|+\sigma, \left( 2l + 2\alpha + 1 \right) \chi_{h}; \htilde \right)}(t).
\label{eqn:cylinder_basis}
\end{equation}
Here $Q_{k}^{(a,b,c;\, p)}$ denotes the degree-$k$ \emph{generalized Jacobi polynomial} orthogonal under the weight $(1-t)^{a} (1+t)^{b} p(t)^{c}$.
Then these basis functions are orthonormal with respect to $\dd \mu(\alpha)$ since
\begin{equation}
\begin{aligned}
\left \langle \Phi_{m,l,k}^{(\alpha,\sigma)}, \Phi_{m',l',k'}^{(\alpha,\sigma)} \right \rangle_{\dd \mu(\alpha)} 
    &= \frac{1}{2 \pi} \int_{0}^{2 \pi} \dd \phi \, e^{-i \left(m - m'\right) \phi} \int_{-1}^{1} \dd \eta \left(1-\eta^{2}\right)^{\alpha} P_{l}^{(\alpha,\alpha)}(\eta) P_{l'}^{(\alpha,\alpha)}(\eta) \\
    &\hspace{-20ex} \times \int_{-1}^{1} \dd t \, \left(1-t\right)^{\left(l + \frac{1}{2} \right) \chi_{o} + \alpha} \left( 1 + t \right)^{\left|m\right| + \sigma} \htilde(t)^{ \left( 2 l + 2 \alpha + 1 \right) \chi_{h} } \\
    &\hspace{-16ex} \times Q_{k}^{\left( \left( l+\frac{1}{2} \right) \chi_{o} + \alpha, \left|m\right|+\sigma, \left( 2l + 2\alpha + 1 \right) \chi_{h}; \htilde \right)}(t)
                          Q_{k'}^{\left( \left( l+\frac{1}{2} \right) \chi_{o} + \alpha, \left|m\right|+\sigma, \left( 2l + 2\alpha + 1 \right) \chi_{h}; \htilde \right)}(t) \\
    &= \delta_{m,m'} \, \delta_{l,l'} \, \delta_{k,k'}.
\end{aligned}
\end{equation}

The full cylinder has a coordinate singularity at $s = 0$.  The basis resolves
this singularity by using the $t$ coordinate, a function of $s^{2}$, in conjunction with the $\left(1+t\right)^{\frac{m}{2}} \propto s^{m}$ prefactor for azimuthal mode $e^{i m \phi}$.
We handle the $s = S_{o}$ coordinate singularity, which only occurs when $\chi_{o} = 1$, by including the $\left( 1-t \right)^{\frac{l}{2}}$ prefactor.  This causes the basis functions of vertical
degree $l$ to decay at the appropriate rate as $t \to 1$ as shown in (\ref{eqn:z_to_the_ell}).

\subsubsection{Differential Operators}
The differential operators in the general cylinder take the form
\begin{equation}
\begin{aligned}
\frac{S_{o}}{2} \mathcal{D}^{+}_{(\alpha,\sigma)} \Phi_{m,l,k}^{(\alpha,\sigma)}
    &=  \phantom{-} \Phi_{m,l,\bigcdot}^{(\alpha+1,\sigma+1)}   \gamma_{l}^{(\alpha)} \bigg[ \left( \mathcal{I}_{h}           \right)^{X_{h}} \, \mathcal{D}(+1,+1,+1) \bigg]_{\bigcdot,k} \\
    &\phantom{=} -  \Phi_{m,l-2,\bigcdot}^{(\alpha+1,\sigma+1)} \delta_{l}^{(\alpha)} \bigg[ \left( \mathcal{I}_{h}^{\dagger} \right)^{X_{h}} \, \mathcal{D}(\delta_{a},+1,-1) \bigg]_{\bigcdot,k} \\
\frac{S_{o}}{2} \mathcal{D}^{-}_{(\alpha,\sigma)} \Phi_{m,l,k}^{(\alpha,\sigma)}
    &=  \phantom{-} \Phi_{m,l,\bigcdot}^{(\alpha+1,\sigma-1)}   \gamma_{l}^{(\alpha)} \bigg[ \left( \mathcal{I}_{h}          \right)^{X_{h}} \, \mathcal{D}(+1,-1,+1) \\
    &\phantom{=} -  \Phi_{m,l-2,\bigcdot}^{(\alpha+1,\sigma-1)} \delta_{l}^{(\alpha)} \bigg[ \left( \mathcal{I}_{h}^{\dagger}\right)^{X_{h}} \, \mathcal{D}(\delta_{a},-1,-1) \bigg]_{\bigcdot,k} \\
\mathcal{D}^{0}_{(\alpha,\sigma)} \Phi_{m,l,k}^{(\alpha,\sigma)}
    &=  \phantom{-} \Phi_{m,l-1,\bigcdot}^{(\alpha+1,\sigma)}   \lambda_{l}^{(\alpha)} \bigg[ \left( \mathcal{I}_{a} \right)^{X_{a}} \bigg]_{\bigcdot,k},
\end{aligned}
\end{equation}
where 
\begin{equation}
\delta_{a} = \begin{cases} +1, \hspace{4ex} \chi_{o} = 1 \\ -1, \hspace{4ex} \chi_{o} = 0 \end{cases},
\hspace{4ex}
X_{a} = \begin{cases} 0, \hspace{4ex} \chi_{o} = 1 \\ 1, \hspace{4ex} \chi_{o} = 0 \end{cases},
\hspace{4ex}
X_{h} = \begin{cases} 1, \hspace{4ex} \chi_{h} = 1 \\ 0, \hspace{4ex} \chi_{h} = \frac{1}{2} \end{cases}.
\end{equation}

Deriving these differential operator coefficients is straightforward but tedious on the surface.
One could differentiate a basis function, thereby computing the grid space differential action of the operator.  They would need to then split up the
operator into its $l \mapsto l + \Delta l$ components using Jacobi polynomial algebra.  Finally they could identify the radial operator for each
$\Delta l$ component, matching it to the analytic expression of the operators as defined in Section~\ref{sec:jacobi_differential_operators}.
This process involves using known Jacobi polynomial algebra rules to simplify conversions like (\ref{eqn:jacobi_conversion_coefficients}).

If we rather follow the required generalized Jacobi parameter increments, ensuring they map appropriately, accounting for
the $s = S_{o}$ coordinate singularity is extremely simple.  Because the $a$ parameter of the 
Jacobi polynomials is $l + \alpha + \frac{1}{2}$ when $\chi_{o} = 1$, the $l \mapsto l-2$ components of the differential operators need an
additional down-conversion of the $a$ parameter.  The generalized Jacobi differential operators naturally achieve this when $\delta_{a} = -1$.

\subsubsection{Coordinate Vector Operators}
The $\vec{s}$ vector multiplication operators are unchanged in the general cylinder.  Multiplication by the axial $\vec{z}$ vector takes the form
\begin{equation}
\vec{z} \, \Phi_{m,l,k}^{(\alpha,0)} = \vec{\hat{e}}_{0} \left ( \beta_{l}^{(\alpha)} \Phi_{m,l+1,\bigcdot}^{(\alpha,0)} \, \bigg[ \left( \mathcal{I}_{a} \right)^{\chi_{o}} \left( \mathcal{I}_{h}^{2} \right)^{\chi_{h}} \bigg]_{\bigcdot,k} + \beta_{l-1}^{(\alpha)}
\Phi_{m,l-1,\bigcdot}^{(\alpha,0)} \, \bigg[ \left( \mathcal{I}_{a}^{\dagger} \right)^{\chi_{o}} \left( \mathcal{I}_{h}^{\dagger \, 2} \right)^{\chi_{h}} \bigg]_{\bigcdot,k} \right).
\end{equation}

\subsubsection{Conversion Operators}
The general conversion operators take the form
\begin{equation}
\begin{aligned}
\mathcal{I}_{(\alpha,\sigma)} \Phi_{m,l,k}^{(\alpha,\sigma)} &\triangleq \Phi_{m,l,k}^{(\alpha,\sigma)} \\
    &= \Phi_{m,l,\bigcdot}^{(\alpha+1,\sigma)}   \, \gamma_{l}^{(\alpha)} \, \bigg[ \mathcal{I}_{a} \,                        \left( \mathcal{I}_{h}^{2}            \right)^{\chi_{h}} \bigg]_{\bigcdot,k}
     - \Phi_{m,l-2,\bigcdot}^{(\alpha+1,\sigma)} \, \delta_{l}^{(\alpha)} \, \bigg[ \left( \mathcal{I}_{a} \right)^{X_{a}} \, \left( \mathcal{I}_{h}^{\dagger \, 2} \right)^{\chi_{h}} \bigg]_{\bigcdot,k}
\end{aligned}
\end{equation}
where
\begin{equation}
X_{a} = \begin{cases} \dagger, \hspace{4ex} \chi_{o} = 1 \\  1, \hspace{4ex} \chi_{o} = 0 \end{cases}.
\end{equation}
That is, for $\chi_{o} = 1$, the $l \mapsto l - 2$ component has the down-conversion operator $\mathcal{I}_{a}^{\dagger}$.
We omit the conversion adjoint here.  To compute them replace $\Delta l$ with $-\Delta l$ and take the adjoint of each radial operator.

\subsection{The General Annulus}
We proceed analogously to the stretched cylinder case, using the appropriate volume element.
Define the weighted measure $\dd \mu(\alpha)$ for $\alpha > -1$ such that
\begin{equation}
\begin{aligned}
\dd \mu (\alpha) &\triangleq \frac{1}{2 \pi} \frac{4}{S_{o}^{2} - S_{i}^{2}} \left( 1-\eta^{2} \right)^{\alpha} \left( 1 - t^{2} \right)^{\alpha} \htilde(t)^{2 \chi_{h} \alpha} \dd V \\
    &= \frac{1}{2 \pi} \left( 1-\eta^{2} \right)^{\alpha} \left( 1 - t \right)^{\frac{1}{2} \chi_{o} + \alpha} \left( 1 + t \right)^{\frac{1}{2} \chi_{i} + \alpha} \htilde(t)^{\left( 2 \alpha + 1 \right) \chi_{h}} \dd \phi \dd \eta \dd t
\end{aligned}
\end{equation}
and induced inner product $\langle \cdot , \cdot \rangle_{\dd \mu(\alpha)}$.

Now we define the gyroscopic basis as follows:
\begin{equation}
\begin{aligned}
\Upsilon_{m,l,k}^{(\alpha,\sigma)} &= e^{i m \phi} \left( S_{i}^{2}(1-t) + S_{o}^{2}(1+t) \right)^{\frac{\left|m\right|+\sigma}{2}} (1-t)^{\frac{l}{2} \chi_{o}} \, (1+t)^{\frac{l}{2} \chi_{i}} \htilde(t)^{l \chi_{h}} P_{l}^{(\alpha,\alpha)}(\eta) \\
    &\hspace{8ex} \times Q_{k}^{\left( \left( l+\frac{1}{2} \right) \chi_{o} + \alpha, \left( l+\frac{1}{2} \right) \chi_{i} + \alpha, \left( 2l + 2\alpha + 1 \right) \chi_{h}, \left|m\right|+\sigma; \htilde, \stilde \right)}(t),
\label{eqn:annulus_basis}
\end{aligned}
\end{equation}
so that 
\begin{equation}
\left \langle \Upsilon_{m,l,k}^{(\alpha,\sigma)}, \Upsilon_{m',l',k'}^{(\alpha,\sigma)} \right \rangle_{\dd \mu(\alpha)} = \delta_{m,m'} \, \delta_{l,l'} \, \delta_{k,k'}.
\end{equation}

Like the stretched cylinder case, the bases naturally conform to all coordinate singularities.  When the domain vanishes at $s = S_{i}$ or $s = S_{o}$ we insert the
$\left( 1 + t \right)^{\frac{l}{2}}$ or $\left( 1 - t \right)^{\frac{l}{2}}$ prefactor, respectively, guaranteeing fast-enough decay for vertical modes of degree $l$
as $s$ approaches the boundary.

\subsubsection{Differential Operators}
The differential operators in the general annulus are
\begin{equation}
\begin{aligned}
\frac{S_{o}^{2} - S_{i}^{2}}{2} \mathcal{D}^{+}_{(\alpha,\sigma)} \Upsilon_{m,l,k}^{(\alpha,\sigma)}
    &=  \phantom{-} \Upsilon_{m,l,\bigcdot}^{(\alpha+1,\sigma+1)}   \, \gamma_{l}^{(\alpha)} \, \bigg[ \left( \mathcal{I}_{h}           \right)^{X_{h}} \, \mathcal{D}(+1,+1,+1,+1) \bigg]_{\bigcdot,k} \\
    &\phantom{=} -  \Upsilon_{m,l-2,\bigcdot}^{(\alpha+1,\sigma+1)} \, \delta_{l}^{(\alpha)} \, \bigg[ \left( \mathcal{I}_{h}^{\dagger} \right)^{X_{h}} \, \mathcal{D}(\delta_{a},\delta_{b},-1,+1) \bigg]_{\bigcdot,k} \\
\frac{S_{o}^{2} - S_{i}^{2}}{2} \mathcal{D}^{-}_{(\alpha,\sigma)} \Upsilon_{m,l,k}^{(\alpha,\sigma)}
    &=  \phantom{-} \Upsilon_{m,l,\bigcdot}^{(\alpha+1,\sigma-1)}   \, \gamma_{l}^{(\alpha)} \, \bigg[ \left( \mathcal{I}_{h}           \right)^{X_{h}} \, \mathcal{D}(+1,+1,+1,-1) \bigg]_{\bigcdot,k} \\
    &\phantom{=} -  \Upsilon_{m,l-2,\bigcdot}^{(\alpha+1,\sigma-1)} \, \delta_{l}^{(\alpha)} \, \bigg[ \left( \mathcal{I}_{h}^{\dagger} \right)^{X_{h}} \, \mathcal{D}(\delta_{a},\delta_{b},-1,-1) \bigg]_{\bigcdot,k} \\
\mathcal{D}^{0}_{(\alpha,\sigma)} \Upsilon_{m,l,k}^{(\alpha,\sigma)}
    &= \phantom{-} \Upsilon_{m,l-1,\bigcdot}^{(\alpha+1,\sigma)}   \, \lambda_{l}^{(\alpha)} \, \bigg[ \left( \mathcal{I}_{a} \right)^{X_{a}} \, \left( \mathcal{I}_{b} \right)^{X_{b}} \bigg]_{\bigcdot,k},
\end{aligned}
\end{equation}
where 
\begin{equation}
\delta_{a,b} = \begin{cases} +1, \hspace{4ex} \chi_{o,i} = 1 \\ -1, \hspace{4ex} \chi_{o,i} = 0 \end{cases},
\hspace{4ex}
X_{a,b} = \begin{cases} 0, \hspace{4ex} \chi_{o,i} = 1 \\ 1, \hspace{4ex} \chi_{o,i} = 0 \end{cases},
\hspace{4ex}
X_{h} = \begin{cases} 1, \hspace{4ex} \chi_{h} = 1 \\ 0, \hspace{4ex} \chi_{h} = \frac{1}{2} \end{cases}.
\end{equation}

Again deriving these operators is natural when we follow the Jacobi parameter increments.  The only modifications compared to the non-singular
geometries are the $\delta_{a,b}$ increments and the $\mathcal{I}_{h}$ powers.  These all can be ``read off'' when we split the operator into
separate $\Delta l$ components and compute the target radial Jacobi parameters.

\subsubsection{Coordinate Vector Operators}
The $\vec{s}$ vector multiplication operators are unchanged in the general annulus.  Multiplication by the axial $\vec{z}$ vector takes the form
\begin{equation}
\begin{aligned}
\vec{z} \, \Phi_{m,l,k}^{(\alpha,0)} &= \vec{\hat{e}}_{0} \bigg( \beta_{l}^{(\alpha)} \Phi_{m,l+1,\bigcdot}^{(\alpha,0)} \, \bigg[ \left( \mathcal{I}_{a} \right)^{\chi_{o}} \left( \mathcal{I}_{b} \right)^{\chi_{i}} \left( \mathcal{I}_{h}^{2} \right)^{\chi_{h}} \bigg]_{\bigcdot,k}  \\
&\hspace{4ex}+ \beta_{l-1}^{(\alpha)} \Phi_{m,l-1,\bigcdot}^{(\alpha,0)} \, \bigg[ \left( \mathcal{I}_{a}^{\dagger} \right)^{\chi_{o}} \left( \mathcal{I}_{b}^{\dagger} \right)^{\chi_{i}} \left( \mathcal{I}_{h}^{\dagger \, 2} \right)^{\chi_{h}} \bigg]_{\bigcdot,k} \bigg).
\end{aligned}
\end{equation}

\subsubsection{Conversion Operators}
The general conversion operators take the form
\begin{equation}
\begin{aligned}
\mathcal{I}_{(\alpha,\sigma)} \Upsilon_{m,l,k}^{(\alpha,\sigma)} &\triangleq \Upsilon_{m,l,k}^{(\alpha,\sigma)} \\
    &= \Upsilon_{m,l,\bigcdot}^{(\alpha+1,\sigma)}   \, \gamma_{l}^{(\alpha)} \, \bigg[ \mathcal{I}_{a} \, \mathcal{I}_{b} \, \left( \mathcal{I}_{h}^{2}            \right)^{\chi_{h}} \bigg]_{\bigcdot,k}
     - \Upsilon_{m,l-2,\bigcdot}^{(\alpha+1,\sigma)} \, \delta_{l}^{(\alpha)} \, \bigg[ \left( \mathcal{I}_{a} \right)^{X_{a}} \, \left( \mathcal{I}_{b} \right)^{X_{b}} \, \left( \mathcal{I}_{h}^{\dagger \, 2} \right)^{\chi_{h}} \bigg]_{\bigcdot,k}
\end{aligned}
\end{equation}
where
\begin{equation}
X_{a,b} = \begin{cases} \dagger, \hspace{4ex} \chi_{o,i} = 1 \\  1, \hspace{4ex} \chi_{o,i} = 0 \end{cases}.
\end{equation}

\subsection{Half Domains are Incompatible with Vanishing Heights}
Placing the domain boundary at $z = 0$ rather than $z = -h(s)$ produces a structural change in the derivative coupling.
This structural change, $\eta \, \deta \mapsto \left( 1 + \zeta \right) \, \dzeta$, induces odd vertical mode coupling in the radial derivatives:
derivatives map $l \mapsto l-1$ in addition to the standard coupling $l \mapsto l$ and $l \mapsto l-2$.
The $l + \frac{1}{2} + \alpha$ Jacobi parameter for sphere-type boundaries and root height functions
does not allow for odd mode coupling.  To see this, note that reaching the $\left( 1 - t \right)^{\frac{l-1}{2}}$ factor present in $l-1$ basis functions
from the  $\left( 1 - t \right)^{\frac{l}{2}}$ factor in $l$ basis functions requires multiplication by $\sqrt{1-t}$, which is not implementable with sparse Jacobi polynomial algebra.
This means the only possible upper half geometries in the $\zeta$ coordinate are standard cylinders and annuli
with polynomial heights.

\end{appendices}

\bibliographystyle{unsrt}
\bibliography{References}

\end{document}